
\input gtmacros
\input rlepsf
\input gtmonout
\volumenumber{1}
\volumeyear{1998}
\volumename{The Epstein birthday schrift}
\input newinsert
\pagenumbers{451}{478}
\received{12 November 1997}
\revised{12 August 1998}
\published{27 October 1998}
\papernumber{22}

\long\def\references{%
\penalty-800\vskip-\lastskip\vskip 15pt plus10pt minus5pt 
{\large\bf References}\ppar      
{\leftskip=25pt\frenchspacing    
\small\parskip=1pt plus1pt       
\def\refkey##1##2\par{\noindent  
\llap{[##1]\stdspace}##2\par}         
\def\,{\thinspace}\thereflist\par}}

\let\ssp\stdspace

\let\Cal\cal

\def\eps{\varepsilon}

\def\R{{\Bbb R}}

\def\slc{{\bf SL}(2,{\C})}

\def\pu2{{\bf PU}(2,1)}

\def\su2{{\bf SU}(2,1)}

\def\po2{{\bf PO}(2,1)}

\def\tr{\mathop{\rm tr}\nolimits}

\def\ax{\mathop{\rm Ax}\nolimits}
\def\Re{\mathop{\rm Re}}
\def\Im{\mathop{\rm Im}}

\def\sech{\mathop{\rm sech}\nolimits}
\def\csch{\mathop{\rm csch}\nolimits}
\def\coth{\mathop{\rm coth}\nolimits}

\mathcode`\:="603A

\def\mod#1{{\left\vert #1 \right\vert}}

\def\gp#1{{\langle #1 \rangle}}

\def\em{\sl}
\def\axis{\mathop{\rm axis}\nolimits}
\def\fix{\mathop{\rm fix}\nolimits}

\def\bound{\partial\,}
\def\PSL{\mathop{\rm PSL}(2,{\Bbb C})}
\def\slc{\mathop{\rm SL}(2,{\Bbb C})}
\def\hat{\widehat}
\def\tilde{\widetilde}

\reflist

\refkey\refAres   
{\bf P Ar\'es},
{\it Coordinates for  Teichm\"uller spaces of $b$-groups with torsion},
Ann. Acad. Sci. Fenn. 20 (1995) 279--300

\refkey\refBus       
{\bf P Buser},
{\it Geometry and Spectra of Compact Riemann Surfaces}, Birkhauser (1992)

\refkey\refEaM    
{\bf C\,J Earle}, {\bf A Marden}, 
{\it Geometric complex coordinates for Teichm\"uller space}, 
in preparation

\refkey\refEpM       
{\bf D\,B\,A Epstein}, {\bf A Marden},
{\it Convex hulls in hyperbolic space, a theorem of Sullivan and 
measured pleated surfaces}, from: ``Analytic and geometric aspects of
hyperbolic space'', D\,B\,A Epstein (editor),
LMS Lecture Notes 111, Cambridge University Press, (1987) 112--253

\refkey\refKSM     
{\bf L Keen}, {\bf C Series},
{\it Pleating coordinates for the Maskit embedding of the %
Teichm\"uller space of punctured tori}, Topology 32 (1993) 719--749

\refkey\refKSC    
{\bf L Keen}, {\bf C Series},
{\it Continuity of convex hull boundaries}, Pacific J. Math. 168
(1995) 183--206

\refkey\refKSQ	
{\bf L Keen}, {\bf C Series},
{\it How to bend pairs of punctured tori}, 
from: ``Lipa's legacy'', J Dodziuk and L
Keen (editors), Contemporary Mathematics 211 (1997) 359--387

\refkey\refKSPQ	
{\bf L Keen}, {\bf C Series},
{\it Pleating invariants for punctured torus groups}, 
{\sl Warwick preprint 10/1998} 

\refkey\refCKB     
{\bf C Kourouniotis},
{\it Deformations of hyperbolic structures}
Math. Proc. Cambridge Phil. Soc. 98 (1985) 247--261

\refkey\refCK       
{\bf C Kourouniotis},
{\it The geometry of bending quasi-Fuchsian groups}, from: ``Discrete
groups and geometry'', W\,J Harvey and C Maclachlan (editors),
LMS Lecture Notes 173, Cambridge University Press, (1992) 148--164

\refkey\refCKC     
{\bf C Kourouniotis},
{\it Complex length coordinates for quasi-Fuchsian groups},
Mathematika 41 (1994) 173--188

\refkey\refKra     
{\bf I Kra},
{\it Horocyclic coordinates for Riemann surfaces and moduli spaces 
I: Teichm\"uller and Riemann spaces of Kleinian groups},
Journal Amer. Math. Soc. 3 (1990) 500--578

\refkey\refKM      
{\bf I Kra}, {\bf B Maskit},
{\it The deformation space of a Kleinian group},
American J. Math. 103 (1980) 1065--1102

\refkey\refBMc     
{\bf B Maskit},
{\it On Klein's combination theorem II},
Trans. AMS 131 (1968) 32--39

\refkey\refBM	
{\bf B Maskit},
{\it Kleinian Groups},
Springer--Verlag, (1987)

\refkey\refMcM	
{\bf C McMullen},
{\it Complex earthquakes and Teichm\"uller theory},
J. Amer. Math. Soc. 11 (1998) 283--320

\refkey\refJN	
{\bf J Nielsen},
{\it Die Isomorphismen der allgemeinen unendlichen Gruppe mit zwei 
Erzeugenden}, Math. Ann. 78 (1918) 385--397

\refkey\refPS      
{\bf J\,R Parker}, {\bf C Series},
{\it Bending formulae for convex hull boundaries},
J. d'Ana\-lyse Math. 67 (1995) 165--198

\refkey\refJouni   
{\bf J Parkkonen},
{\it Geometric complex analytic coordinates for deformation spaces 
of Koebe groups},	
Ann. Acad. Sci. Fenn. Math. Diss. 102 (1995) 1--50    

\refkey\refSPT    
{\bf S\,P Tan},
{\it Complex Fenchel-Nielsen coordinates for quasi-Fuchsian 
structures}, International J. Math. 5 (1994) 239--251

\refkey\refWW      
{\bf P\,L Waterman}, {\bf S\,A Wolpert},
{\it Earthquakes and tessellations of Teichm\"uller space},
Trans. AMS 278 (1983) 157--167

\refkey\refSW	
{\bf S\,A Wolpert},
{\it The Fenchel--Nielsen deformation},
Annals of Math. 115 (1982) 501--528

\refkey\refSWS	
{\bf S\,A Wolpert},
{\it On the symplectic geometry of deformations of a hyperbolic 
surface}, Annals of Math. 117 (1983) 207--234

\refkey\refDW 
{\bf D\,J Wright},
{\it The shape of the boundary of Maskit's embedding of the 
Teichm\"uller space of once punctured tori}, preprint

\endreflist

\title{Coordinates for Quasi-Fuchsian\\Punctured Torus Space}

\author{John R Parker\\Jouni Parkkonen}

\abstract
We consider complex Fenchel--Nielsen coordinates on the\break
quasi-Fuchsian space of punctured tori. These coordinates arise from a
generalisation of Kra's plumbing construction and are related to
earthquakes on Teichm\"uller space. They also allow us to interpolate
between two coordinate systems on Teichm\"uller space, namely the
classical Fuchsian space with Fenchel--Nielsen coordinates and the
Maskit embedding.  We also show how they relate to the pleating
coordinates of Keen and Series.
\endabstract

\asciiabstract{%
We consider complex Fenchel-Nielsen coordinates on the quasi-Fuchsian
space of punctured tori. These coordinates arise from a generalisation
of Kra's plumbing construction and are related to earthquakes on
Teichmueller space. They also allow us to interpolate between two
coordinate systems on Teicmueller space, namely the classical Fuchsian
space with Fenchel-Nielsen coordinates and the Maskit embedding.  We
also show how they relate to the pleating coordinates of Keen and
Series.}

\address{Department of Mathematical Sciences, University of Durham\\
Durham DH1 3LE, UK\\\smallskip\\
Department of Mathematics, University of Jyv\"askyl\"a\\
P.O.Box 35, FIN-40351 Jyv\"askyl\"a, Finland}

\asciiaddress{Department of Mathematical Sciences, University of Durham\\
Durham DH1 3LE, England,
Department of Mathematics, University of Jyvaskyla\\
P.O.Box 35, FIN-40351 Jyvaskyla, Finland}


\email{J.R.Parker@durham.ac.uk, parkkone@math.jyu.fi}

\maketitle
\sectionnumber=-1
\section{Introduction}

In this note we study the holomorphic extension of the classical
Fenchel--Nielsen coordinates of the Teichm\"uller space of
once-punctured tori to the quasi-con\-formal deformation space of a
Fuchsian group representing two punctured tori, quasi-Fuchsian
punctured torus space. A {\em punctured torus group} $G=\langle
S,T\rangle$ is a discrete, marked, free subgroup of $\PSL$ with two
generators whose commutator $K=T^{-1}S^{-1}TS$ is parabolic.  This
group acts naturally on the Riemann sphere by conformal
transformations.  The limit set $\Lambda(G)$ consists of all
accumulation points of this action and is the smallest nonempty closed
$G$--invariant subset of the Riemann sphere. Its complement is called
the ordinary set $\Omega(G)$. The group $G$ is called {\em
quasi-Fuchsian} if its ordinary set $\Omega(G)$ consists of two simply
connected components or equivalently if its limit set $\Lambda(G)$ is
a topological circle. The space of all quasi-Fuchsian punctured torus
groups up to conjugation within $\PSL$ is called {\em quasi-Fuchsian
punctured torus space} and will be denoted by ${\Cal Q}$. The subset
of ${\Cal Q}$ consisting of groups whose limit set is a round circle
is the space of all Fuchsian punctured torus groups.  We call this
{\em Fuchsian punctured torus space} and we will denoted it by ${\Cal
F}$. It is a copy of the Teichm\"uller space of the punctured torus.

Our approach to quasi-Fuchsian punctured torus groups is a combination of
the classical Fenchel--Nielsen construction of Fuchsian groups and the
gluing construction used by Kra in [\refKra] for terminal $b$-groups. This
is rather natural as Fuchsian groups form a real subspace inside the
space of quasi-Fuchsian groups, and terminal $b$-groups form part of the
boundary of the same space. We start with a Fuchsian group $F$ of the
second kind such that $X_0$, the quotient of the hyperbolic plane by $F$,
is a sphere 
with a puncture and two infinite area ends with boundary geodesics of equal
lengths. We then extend the group by adding a M\"obius transformation
that glues together the infinite area ends of the quotient to make a 
punctured torus. If the resulting group $G$ is Fuchsian, this is the
Fenchel--Nielsen construction. The construction is carried out in Section~1
and the Fenchel--Nielsen parameter is connected with the gluing
parameter in Proposition~3.2. We can also regard $F$ and $G$ as acting on the
Riemann sphere and we allow the Fenchel--Nielsen parameters to be complex.
For other allowed values of the gluing parameter
the resulting group $G$ is a quasi-Fuchsian group bent along the geodesic in
${\Bbb H}^3$ corresponding to the boundary geodesics of $X_0$. The
analysis of this bending, the associated shear, and their use for
parametrising the deformation space of quasi-Fuchsian groups from different
points of view is the main
goal of the second half of the paper. We show that the resulting complexified 
Fenchel--Nielsen twist parameter can be interpreted as a complex shear as 
introduced by Parker and Series in [\refPS] and that it has another natural 
interpretation as a $zw=t$ plumbing parameter as in Kra [\refKra] 
The relationship between the various points of view is often easy at a
conceptual level but can be hard to make explicit. In this paper we aim to
make these connections as explicit as possible. Part of this involves writing 
down generators for punctured torus groups as matrices depending on 
parameters. This is useful for making explicit computations which we illustrate
by drawing pictures of various slices through ${\cal Q}$.

One of the main themes of this paper will be a partial description of
Keen--Series pleating invariants in terms of complex Fenchel--Nielsen
parameters. For completeness we now give a brief account of pleating 
invariants [\refKSM, \refKSPQ]. Unlike complex 
Fenchel--Nielsen coordinates these are not holomorphic coordinates but they
do reflect the geometrical structure of the associated 3--manifold as well
as the limit set of $G$. In particular, they may be used to determine the
shape of the embedding of ${\cal Q}$ into ${\Bbb C}^2$ given by complex
Fenchel--Nielsen coordinates. We will illustrate this with pictures of various 
slices through this embedding. Let $G$ be a punctured torus group that is 
quasi-Fuchsian but not Fuchsian. We call such a group 
{\em strictly quasi-Fuchsian}. Consider $C(G)$, the the hyperbolic convex hull 
in ${\Bbb H}^3$ of the limit set of $G$ (sometimes called the Nielsen region 
for $G$). 
This is a $G$--invariant, simply connected, convex subset of ${\Bbb H}^3$.
Thus, its quotient $C(G)/G$ is a convex 3--manifold with boundary,
whose fundamental group is $G$. In other words $C(G)/G$ is topologically, the 
product of a closed interval with a punctured torus. Each boundary component 
is topologically a punctured torus and naturally inherits a hyperbolic
structure from the three manifold (this structure is different from the
obvious hyperbolic structure on the corresponding component of 
$\Omega(G)/G$). This hyperbolic structure makes the boundary component
into a pleated surface in the sense of Thurston. That is, it consists of
totally geodesic flat pieces joined along a geodesic lamination, called
the {\em pleating locus}, and which carries a natural transverse measure, the
{\em bending measure}. The length $l_\mu$ of a measured lamination $\mu$
on a surface with a given hyperbolic structure, is the total mass 
on this surface of the measure given by the product of hyperbolic length
along the leaves of $\mu$ with the transverse measure $\mu$. 
For the punctured torus it is well known that measured geodesic laminations
are projectively parametrised by the extended real line. If the support of 
the lamination is drawn on the square flat torus then this parameter is just 
the gradient. From this we see that the possible types of support that this 
lamination that can have fall into two categories. First, simple closed 
curves, sometimes called
rational laminations because of their parametrisation by rational
slopes on a square torus. The transverse measure is just the $\delta$--measure
on these curves. Secondly, laminations whose leaves are
unbounded geodesic arcs and which correspond to ``infinite words'' in $G$.
We refer to these as {\em infinite laminations}. They correspond to curves of
irrational slope on a square torus and so are sometimes referred to as
irrational laminations. The measure they carry is called bending measure.
We remark that the pleating locus cannot be the same on both components of
the convex hull boundary. This is an important observation.
Most of the time in this paper, we will be concerned with
the case where the pleating locus on one component of $\partial C(G)/G$ is a 
simple closed geodesic. In this case,
there will be a constant angle across this geodesic between the two
adjacent flat pieces. In this case, the lamination length is just the
length of the geodesic in the hyperbolic structure on the convex hull
boundary. Keen and Series show in [\refKSPQ] that a marked punctured
torus group is determined by its pleating invariants, namely the projective
classes $(\mu,l_\mu)$, $(\nu, l_\nu)$ where the supports of $\mu$ and $\nu$ 
are the pleating loci on the two components of $\partial C(G)$ and $l_\mu$, 
$l_\nu$ are their lamination lengths.
 
Suppose that the pleating loci on both components of $\partial C(G)$ are
simple closed curves $\gamma$, $\delta$. The corresponding group elements
necessarily have real trace (though this is not a sufficient condition).
The collection of all groups in ${\cal Q}$ for which $\gamma$, $\delta$ are
the pleating loci is called the {\em (rational) pleating plane} 
${\cal P}_{\gamma,\delta}$. This is a two dimensional non-singular 
subset of ${\cal Q}$ and is parametrised by the lengths of the geodesics
$\gamma$ and $\delta$ (which in this case are the lamination lengths),
see Theorem~2 of [\refKSPQ]. Keen and Series also define pleating planes
for the cases where one or both of the pleating loci are infinite laminations.
We will only make passing reference to such pleating planes.

We have been greatly helped by conversations with Linda Keen and Caroline
Series. We would like to thank them for their help. We would also like to
thank the referee for her/his comments which have improved the paper.
The second author was supported by the Academy of Finland and by the
foundation Magnus Ernroothin S\"a\"ati\"o of the Finnish Society of Sciences
and Letters.
Figures 4.1, 5.1 and 6.1 were drawn using a computer program developed by
David Wright. The second author would like to thank him for his help in
installing and using the program.
Both authors would like to thank the Centre Emile Borel at the Institut Henri
Poincar\'e for their hospitality.

\section{Real Fenchel--Nielsen coordinates}

In this section we show how to write down generators for Fuchsian punctured
torus groups in terms of Fenchel--Nielsen coordinates. This section gives 
a foundation for the subsequent sections: In order to obtain complex
Fenchel--Nielsen coordinates we simply keep the same normal form for the 
generators but make the parameters complex. The material in this section is
quite standard, for a more complete discussion of  Fenchel--Nielsen
coordinates see Buser [\refBus]. 

Let $X$ be a punctured torus and
$\gamma\subset X$ a simple closed
geodesic. Then $X_0=X\setminus\gamma$ is a hyperbolic surface of genus $0$
with
one puncture and two geodesic boundary components of equal length, say $l$.
$X_0$ can be realised as a quotient $X_0=N(G_0)/G_0$, where 
$G_0$
is a
Fuchsian group of
the second kind generated by two hyperbolic transformations with multiplier
$\lambda=l/2\in\R_+$:
$$
S=\pmatrix{ \cosh(\lambda) & \cosh(\lambda)+1 \cr \cosh(\lambda)-1 &
\cosh(\lambda)}\quad{\rm and}
\quad
S'=\pmatrix{ \cosh(\lambda) & \cosh(\lambda)-1 \cr \cosh(\lambda)+1 &
\cosh(\lambda) }, \eqno{(1.1)}
$$
and $N(G_0)$ is the {\em Nielsen region} of $G_0$, that is, the
hyperbolic convex hull in ${\Bbb H}$ of the limit set of $G_0$. For later
reference we record that the fixed points of these transformations are 
$\fix S=\pm\coth(\lambda/2)$ and $\fix S'=\pm\tanh(\lambda/2)$. 
The transformations $S$ and $S'$ correspond to the boundary geodesics of
$X_0$ and their product $K={S'}^{-1}S$ corresponds to the puncture.
In other words
$$
K={S'}^{-1}S=\pmatrix{ -1+2\cosh(\lambda) & 2\cosh(\lambda) \cr
-2\cosh(\lambda) & -1-2\cosh(\lambda) }\eqno{(1.2)}
$$
is a parabolic transformation fixing $-1$.

\midinsert
\centerline{\relabelbox\small
\epsfxsize 10 truecm\epsfbox{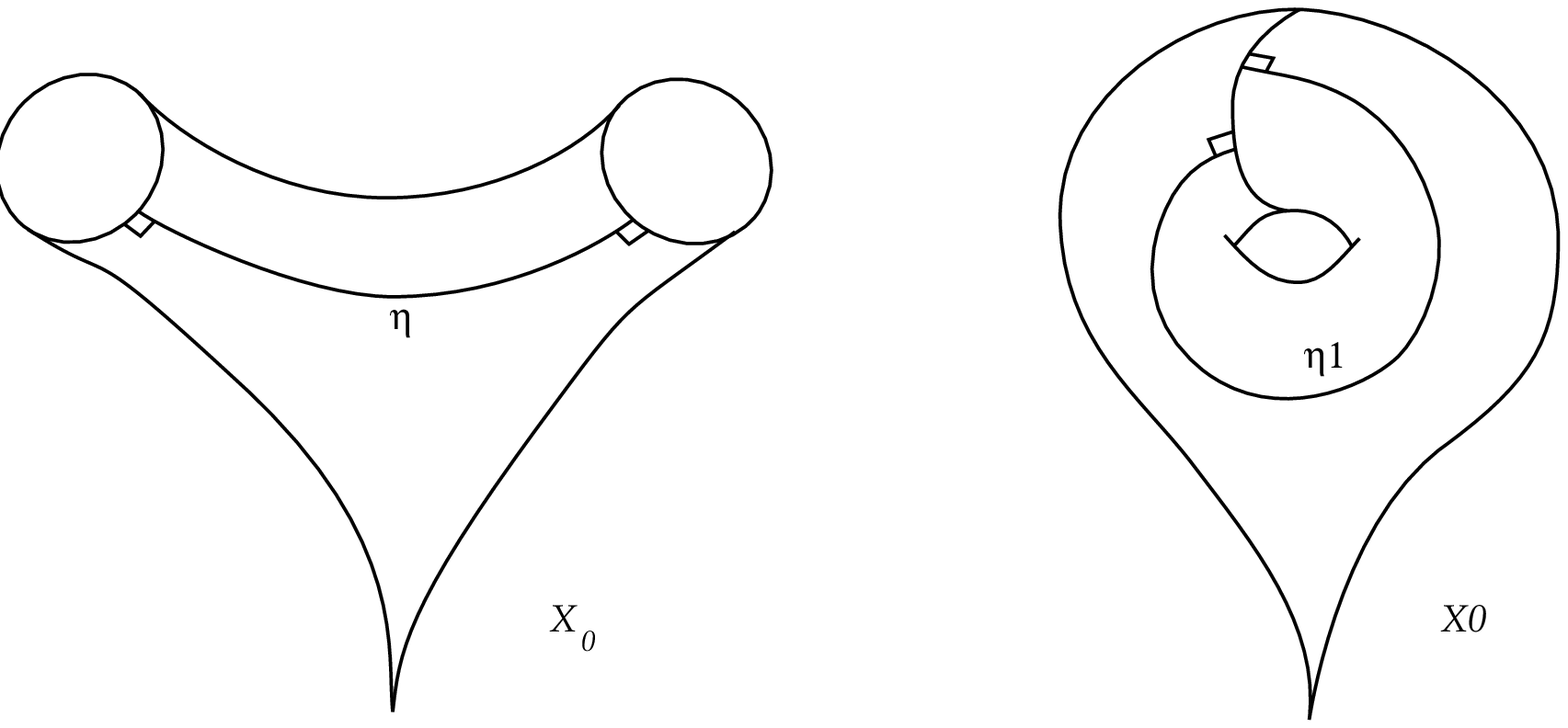}
\relabel {h}{$\eta$}
\relabel {h1}{$\eta$}
\relabel {X}{$X_0$}
\relabel {X0}{$X$}
\endrelabelbox}
\centerline{\small Figure~1.1\ssp The Fenchel--Nielsen construction}
\endinsert

The original surface $X$ can be reconstructed by gluing together the geodesic
boundary components of $X_0$. The gluing can be realised by adding to the 
group a hyperbolic M\"obius transformation $T$ that preserves ${\Bbb H}^2$.
We form a new Fuchsian group, an HNN extension of $G_0$:
$$
G=\langle G_0,T\rangle=\left(G_0\right)*_{\langle T\rangle}.
$$
The transformation $T$ is required to conjugate the cyclic subgroups 
$\langle S\rangle$ and $\langle S'\rangle$ in a manner compatible with the 
gluing operation:
$$
T^{-1}ST=S'.
$$
This condition fixes $T$ up to one free parameter $\tau\in\R$, and $T$ can
be written in the form
$$
T=\pmatrix{ \cosh(\tau/2)\coth(\lambda/2) & -\sinh(\tau/2) \cr
-\sinh(\tau/2) & \cosh(\tau/2)\tanh(\lambda/2) }.\eqno{(1.3)}
$$
\midinsert
\centerline{
\relabelbox
\epsfxsize 12.5 truecm\epsfbox{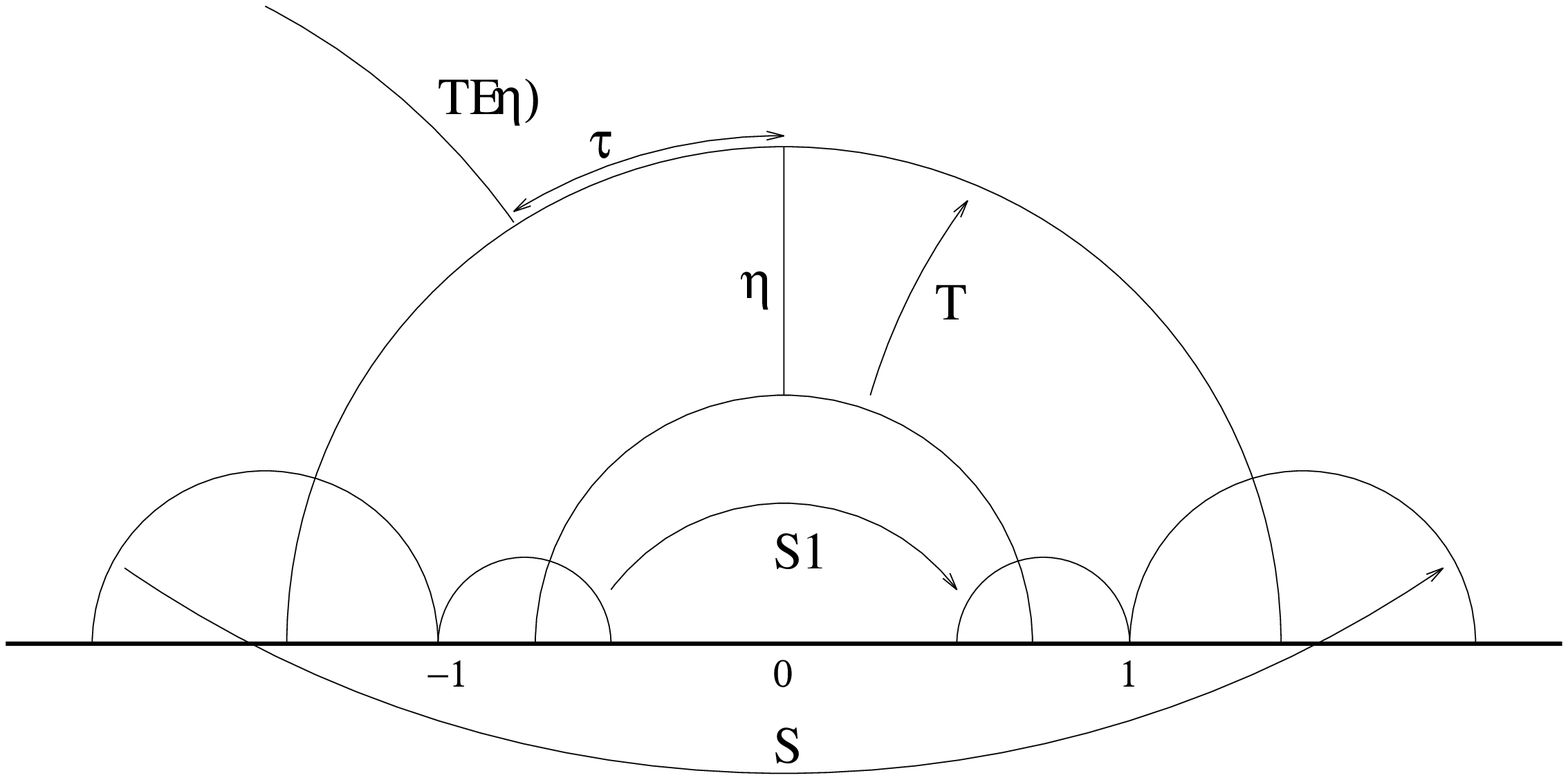}
\relabel {T}{\small$T$}
\relabel {TE}{\small$T(\eta)$}
\relabel {t}{\small$\tau$}
\relabel {h}{\small$\eta$}
\relabel {S}{\small$S$}
\relabel {S1}{\small$S'$}
\relabel {1}{$\scriptstyle1$}
\relabel {-1}{$\scriptstyle-1$}
\relabel {0}{$\scriptstyle0$}
\endrelabelbox}
\vglue 0.05truein
\centerline{\small Figure~1.2\ssp The fundamental domain}
\endinsert
We recover the original (marked) surface with the correct geometry for
exactly  one parameter $\tau_0\in\Bbb R$.
However, the group $G$ is a Fuchsian group for any real $\tau$, and the
parameter has a geometric interpretation:
There is a unique simple geodesic arc $\eta$ on $X_0$ perpendicular to both
geodesic boundary curves. A distinguished lift of this arc to the universal
covering ${\Bbb H}^2$ is the segment of the positive imaginary axis
connecting $i\tanh(\lambda/2)\in\axis(S')$ and
$i\coth(\lambda/2)\in\axis(S)$.
Now $T$ maps $i\tanh(\lambda/2)$ to a point on the axis of $S$,
namely 
$$
T(i\tanh(\lambda/2))=
i\coth(\lambda/2)\bigl(\sech(\tau)+i\tanh(\tau)\bigr).
$$
The (signed) hyperbolic distance of this point from
$i\coth(\lambda/2)$ is exactly $\tau$, 
the sign of $\tau$ is chosen to be positive if moving from $i\coth(\lambda/2)$
to $T\bigl(i\tanh(\lambda/2)\bigr)$ takes one in a positive (anti-clockwise)
direction around the circle of radius $\coth(\lambda/2)$.
The map $G\longmapsto(\lambda,\tau)$ is the {\em Fenchel--Nielsen coordinate} 
of the Teichm\"uller space of punctured tori. It defines a global 
real analytic parametrisation and identifies ${\Cal F}$ with 
${\Bbb R}_+\times\Bbb R$ (see Buser [\refBus]).
Fenchel--Nielsen coordinates depend on the choice of an ordered pair of 
(homotopy classes of) simple closed curves on the punctured torus intersecting 
exactly once, that is a marking. We obtain different coordinates for 
different choices of marking. These choices are related by elements of the 
modular group. We investigate this in more detail in the next section.
In [\refWW] Waterman and Wolpert give computer pictures for the action of
the modular group on Fenchel--Nielsen coordinates. They also give
pictures of this action in another set of coordinates which can be easily 
derived from traces of generating triples. 

Varying $\tau$ and keeping $\lambda$ fixed is the Fenchel--Nielsen deformation
considered by Wolpert in [\refSW] and [\refSWS].

\section{Complex Fenchel--Nielsen coordinates}

The Teichm\"uller space of punctured tori seen as the space of Fuchsian groups
representing a punctured torus, ${\Cal F}$, is a natural subspace of
the corresponding quasi-Fuchsian space, $\Cal Q$. 
Kourouniotis [\refCKC] and Tan [\refSPT] showed that, for compact surfaces, 
the Fenchel--Nielsen coordinates can be complexified to give a global 
parametrisation of quasi-Fuchsian space.
With this in mind we now suppose that $\lambda$ and $\tau$ are complex.
That is $(\lambda,\tau)\in{\Bbb C}_+\times{\Bbb C}$ where
${\Bbb C}_+$ denotes those complex numbers with positive real part.
With such $\lambda$ and $\tau$ we consider groups generated by $S$ and $T$ 
with the normal forms (1.1) and (1.3).
This means that $S$ and $T$ are now in $\PSL$ rather than in 
${\rm PSL}(2,{\Bbb R})$.
The group $\gp{S,T}$ is not quasi-Fuchsian for
all $(\lambda,\tau)\in{\Bbb C}_+\times{\Bbb C}$ but the {\em complex 
Fenchel--Nielsen coordinates} $(\lambda,\tau)$ do give global coordinates on
${\Cal Q}$. We present a short proof of this fact
using the stratification method developed by Kra and Maskit in [\refKM].

\proclaim{Proposition~2.1}The map $h\co {\Cal Q}\to{\Bbb C}^{2}$ given by 
$h(G)=(\cosh^2(\lambda),e^\tau)$ is 
a global complex analytic coordinate map on $\Cal Q$.\endproc

\prf Let $G=\langle A,B\rangle$ be a quasi-Fuchsian group of type $(1,1)$
generated by two loxodromic transformations $A$ and $B$. Assume that the
group is normalised so that $0$ is the repelling fixed point, and $\infty$ is
the attracting fixed point of $A$, and that $B(0)=1$. Let $x_1=B(\infty)$, and
$x_2=B(1)$. Note that $x_1,x_2\in\Lambda(G)$.

We claim that $G$ is determined by giving $x_1$ and $x_2$:
Clearly $B$ is determined, as we know how it maps three points. Also,
from the normalisation we know that 
$$
A=\pmatrix {a &0\cr 0&1/a},\qquad
B=\pmatrix { x_1(x_2-1) & x_1-x_2 \cr x_2-1 & x_1-x_2 },
$$
where $a\in\Bbb C$, $|a|>1$. Now
$$
\tr [A,B]={{2 a^2 x_1 - 1 - a^4}\over{a^2(x_1-1)}}.
$$
As  $[A,B]$ is assumed to be a parabolic, solving for $a^2$ in the equation
$\tr [A,B]=-2$ gives $a^2=2x_1 - 1\pm 2\sqrt{x_1(x_1-1)} $. 
Only one of these solutions satisfies $|a|>1$. This fixes $A$.
(The choice of the branch of the square root $a=\sqrt{a^2}$ does 
not affect $A$.)

Let us normalise the group $G=\langle S',T\rangle$ of Section 1 as above:
We conjugate $G$ with a transformation (here written as an element of
$\mathop{\rm PGL}(2,\Bbb C)$)
$$
R=
\pmatrix{
\cosh(\lambda)/(1-\cosh(\lambda)) & -\coth(\lambda) \cr
1/(1-\cosh(\lambda)) & \csch(\lambda)}.
$$
This gives
$$
S_0=RS'R^{-1}=\pmatrix{e^\lambda &0\cr 0& e^{-\lambda}},
$$
where we can assume $|e^\lambda|>1$, and 
$$
T_0=RTR^{-1}=\pmatrix{
\coth(\lambda)e^{-\tau/2} & \coth(\lambda)e^{\tau/2} \cr
\csch(\lambda)\sech(\lambda)e^{-\tau/2} & \coth(\lambda)e^{\tau/2}}.
$$
Now 
$$
x_1=\cosh^2(\lambda),\quad\quad x_2={1+e^\tau \over
\sech^2(\lambda)+e^\tau}.
$$
\vglue-.3truein\endprf

\rk{Remark~2.2}The choice $|e^\lambda|>1$ implies $\lambda\in\Bbb C_+$.
Unlike real Fenchel--Nielsen coordinates, there is no simple description
of which pairs $(\lambda,\tau)\in{\Bbb C}_+\times{\Bbb C}$ are in $h(\Cal Q)$,
the image of quasi-Fuchsian space under the coordinate map. Using the pleating
invariants of Keen and Series [\refKSPQ] one can determine how $h(\Cal Q)$ lies
inside ${\Bbb C}^2$. In this paper we carry out part of this construction and
illustrate our results by drawing slices through ${\Cal Q}$ in Figure~5.1.

We now use the fact that $(\cosh^2(\lambda),e^\tau)$ give global coordinates
to show that $(\lambda,\tau)$ give global coordinates on quasi-Fuchsian
space. Let 
$$
\tilde{\Cal{FN}}=\left\{(\lambda,\tau)\in{\Bbb C}^2 :
(\cosh^2(\lambda),e^\tau)\in h(\Cal Q)\right\},
$$
where $h$ is the map of Proposition~2.1. We denote by ${\Cal{FN}}$ the 
component of $\tilde{\Cal{FN}}$ containing $\Bbb R_+\times \Bbb R$. 
Our proof that $(\lambda,\tau)$ give global coordinates involves showing that 
there are no paths in $\tilde{\Cal FN}$ between two places where the 
parameters are different but the groups are the same. 

\proclaim{Proposition~2.3} Let 
$\gamma\co [0,1]\longrightarrow {\Bbb C_+}\times{\Bbb C}$ denote any path from 
$\gamma(0)=(\lambda_0,\tau_0)$ to 
$\gamma(1)=(\lambda_0+m\pi i,\tau_0+2n\pi i)$ for any 
$(\lambda_0,\tau_0)\in{\tilde{\Cal{FN}}}$ and integers 
$m$ and $n$ not both zero. 
Then $\gamma([0,1])$ is not contained in $\tilde{\Cal{FN}}$.\endproc

\prf
We begin with the case $m=1$ and $n=0$.

Using the normalisation of Proposition~2.1 we have
$T_0(\lambda_0,\tau_0)=T_0(\lambda_0+\pi i,\tau_0)$. Also notice that
$S_0(\lambda_0,\tau_0)$ and $S_0(\lambda_0+\pi i,\tau_0)$ are the same in
$\PSL$ but differ by $-I$ in $\slc$. They correspond to the two choices
of square root for $a^2$ in Proposition~2.1. Thus moving
along $\gamma$ from $(\lambda_0,\tau_0)$ to $(\lambda_0+\pi i,\tau_0)$ 
adds $i\pi$ to the multiplier of $S_0$. For more
details of the relationship between multipliers and the different lifts of 
M\"obius transformations in $\PSL$ to matrices in $\slc$ see the discussion
in Section 1 of [\refPS]. 
Let $\Pi_1$ be any hyperplane in ${\Bbb H}^3$ orthogonal
to the axis of $S_0$ and let $\Pi_2=S_0(\Pi_1)$ be its image under $S_0$. 
Because going along $\gamma$ from $(\lambda_0,\tau_0)$ to 
$(\lambda_0+\pi i,\tau_0)$ changes the multiplier of $S_0$ by $\pi i$ then
also $\Pi_2$ is rotated by $2\pi$ with respect to $\Pi_1$. We can think 
of going along $\gamma$ as being the same as doing a Dehn twist of the 
annulus between $\partial\Pi_1$ and $\partial\Pi_2$ in $\widehat{\Bbb C}$.

Specifically we may decompose $S_0$ into a product of half turns (that is
elliptic involutions in $\PSL$ of order 2) as follows:
$$
S_0=\iota_1\iota_2=\pmatrix{0 & e^\lambda \cr -e^{-\lambda} & 0}
\pmatrix{0 & -1 \cr 1 & 0}.
$$
The geodesic fixed by $\iota_1(\lambda,\tau)$ has end points $\pm ie^\lambda$.
Replacing $(\lambda_0,\tau_0)$ by $(\lambda_0+\pi i,\tau_0)$ 
interchanges these end points. Equivalently this reverses the orientation of
the geodesic. Therefore if $\Pi_1$ is hyperplane orthogonal 
to the axis of $S_0$ and containing the geodesic with end points $\pm i$ 
(that is the axis of $\iota_2$) it is clear that its image under the $\iota_1$ 
is rotated by $2\pi$ when we replace $\lambda_0$ by $\lambda_0+\pi i$.

Let $\xi_1$ be any point of $\partial\Pi_1\cap\Omega$ and $\xi_2=S_0(\xi_1)$
be its image under $S_0$. Let $\alpha$ be any path in $\Omega$ joining 
$\xi_1$ and $\xi_2$. Now consider the homotopy $H$ given by following $\alpha$
while $(\lambda,\tau)$ varies along $\gamma$. Denote the image of $\alpha$
at time $t$ by $\alpha_t$.

If the whole of $\gamma$ were in ${\Cal Q}$ then the homotopy
$H$ would induce an isotopy from $\Omega\bigl(G(\lambda_0,\tau_0)\bigr)$
to $\Omega\bigl(G(\lambda_0,\tau_0+2\pi i)\bigr)$. At each stage $S_0$ is 
loxodromic so $\Pi_1$ and $\Pi_2$ are disjoint and $\alpha_t$ consists of
more than one point. Now $\alpha_0$ and $\alpha_1$ are both paths in
$\Omega\bigl(G(\lambda_0,\tau_0)\bigr)=
\Omega\bigl(G(\lambda_0+\pi i,\tau_0)\bigr)$ joining $\xi_1$ and $\xi_2$.
It is clear from the earlier discussion that the path 
$\alpha_1{\alpha_0}^{-1}$ formed by going along $\alpha_1$ and then backwards 
along $\alpha_0$ winds once around the (closed) annulus between 
$\partial\Pi_1$ and $\partial\Pi_2$. This it separates the fixed points of 
$S_0$. This contradicts the fact that the limit set $\Lambda$ is connected. 

We can adapt this proof to cover the case where
$\lambda_0$ is sent to $\lambda_0+m\pi i$ for some non-zero integer $m$.
This is done by observing that the path $\alpha_1{\alpha_0}^{-1}$ now 
winds $m$ times around the annulus between $\partial\Pi_1$ and
$\partial\Pi_2$. Moreover this argument does not use the value of
$\tau$ at each end of the path. It merely uses the fact that 
$T_0(\lambda_0,\tau_0)=T_0(\lambda_1,\tau_1)$ and so we may take
$\tau_1=\tau_0+2n\pi i$ without changing anything. 

Thus we have proved the result when $m$ and $n$ are any integers with $m$ not
zero. It remains to prove the result when $m=0$ and $n$ is an integer other 
than
zero. We do this as follows. Observe that, with the normalisation of (1.1) and 
(1.3), $S(\lambda_0,\tau_0)=S(\lambda_0,\tau_0+2\pi i)$ but
$T(\lambda_0,\tau_0+2\pi i)$ and $T(\lambda_0,\tau_0)$ give distinct lifts
in $\slc$. As before we decompose $T$ into a product of half turns as follows:
$$
T=\iota_1\iota_2=\pmatrix{\sinh(\tau/2) & \cosh(\tau/2)\coth(\lambda/2) \cr
-\cosh(\tau/2)\tanh(\lambda/2) & -\sinh(\tau/2)}
\pmatrix{0 & -1 \cr 1 & 0}.
$$
The geodesic fixed by $\iota_1(\lambda,\tau)$ has end points
$$
{-\sinh(\tau/2)\pm i\over \cosh(\tau/2)\tanh(\lambda/2)}.
$$
Replacing $(\lambda_0,\tau_0)$ by $(\lambda_0,\tau_0+2\pi i)$ 
interchanges these end points. The rest of the argument follows as before.
\endprf

The next two results are direct consequences of Propositions~2.1
and~2.3.

\proclaim{Corollary}The functions $\cosh^2(\lambda)$ and $e^\tau$
have well defined inverses in ${h(\Cal Q)}$ and so we can regard 
$(\lambda,\tau)$ is a global coordinate system for quasi-Fuchsian space. 
\endproc
\proclaim{Corollary}The pair $\bigl(\cosh(\lambda),\sinh(\tau/2)\bigr)$ 
give global coordinates for quasi-Fuchsian space. In particular, the 
points where $\sinh(\lambda)=0$ or $\cosh(\tau/2)=0$ are not in ${\Cal{FN}}$. 
\endproc

\prf 
The first part follows from the previous corollary.
We give a simple justification for the last statement. If 
$\sinh(\lambda)=0$ then $\cosh(\lambda)=\pm 1$ and $S$ is parabolic.
Similarly if $\cosh(\tau/2)=0$ then $T$ is elliptic or else $\coth(\lambda)$
is infinite and $S$ is parabolic as before. \endprf

Complex Fenchel--Nielsen coordinates depend on the choice of a marking for
the punctured torus, that is an ordered pair of generators for $S$. It is
intuitively clear that changing this marking gives a biholomorphic change of 
the coordinates $\bigl(\cosh(\lambda),\sinh(\tau/2)\bigr)$. We now make this
explicit.

\proclaim{Proposition~2.4}Let $(S_0,T_0)$ and $(S_1,T_1)$ be any two 
generating 
pairs for a punctured torus group $G$. Let $(\lambda_0,\tau_0)$ and 
$(\lambda_1,\tau_1)$ be the corresponding complex Fenchel--Nielsen coordinates
on ${\Cal Q}$. Then the map 
$$\bigl(\cosh(\lambda_0),\sinh(\tau_0/2)\bigr)\longmapsto
\bigl(\cosh(\lambda_1),\sinh(\tau_1/2)\bigr)$$
is a biholomorphic homeomorphism of ${\Cal Q}$ to itself.\endproc

\prf
A classical result of Nielsen [\refJN] states that we can obtain the pair 
$(S_1,T_1)$ from $(S_0,T_0)$ by a sequence of elementary Nielsen-moves on the
generators. As one of our aims is to make things explicit, we list these
Nielsen moves and write down the effect that they have on the 
coordinates $\bigl(\cosh(\lambda),\sinh(\tau/2)\bigr)$. From this, it is clear 
that these changes of coordinate are holomorphic.

First, suppose that $(S',T')=(S,S^{\pm 1}T)$. Then
$$
\cosh(\lambda')=\cosh(\lambda),\quad 
\sinh(\tau'/2)=\sinh(\tau/2)\cosh(\lambda)\mp\cosh(\tau/2)\sinh(\lambda).
$$
Secondly, suppose that $(S',T')=(S,T^{-1})$. Then
$$
\cosh(\lambda')=\cosh(\lambda),\quad 
\sinh(\tau'/2)=-\sinh(\tau).
$$
Finally, suppose that $(S',T')=(T,S)$
$$
\cosh(\lambda')={\cosh(\lambda)\cosh(\tau/2)\over\sinh(\lambda)},\quad
\sinh(\tau'/2)={-\sinh(\tau/2)\sinh(\lambda)\over\cosh(\tau/2)}.
$$
\vglue-.3truein\endprf

\section{Plumbing and earthquakes}

In this section we show how the Fenchel--Nielsen construction is related to
two standard constructions in Teichm\"uller theory, namely the $zw=t$
plumbing construction and to quake-bends. In particular, the Fenchel--Nielsen
twist parameter is a special case of the quake-bend parameter and we 
show how to express the plumbing parameter in terms of Fenchel--Nielsen 
parameters.

Consider Teichm\"uller space of the punctured torus $\Cal F$ with 
Fenchel--Nielsen coordinates as in Section 1. The motion 
through Teichm\"uller space obtained by fixing the length parameter 
$\lambda$ but varying the shear $\tau$ is the Fenchel--Nielsen deformation
(see [\refSW]) which is the simplest example of an earthquake
(see Waterman and Wolpert [\refWW] and McMullen [\refMcM] for some other 
earthquakes). 
One may think of this as cutting along $\ax(S)$ twisting and then regluing.

If we reglue so that along $\ax(S)$ the two sides make a constant angle then 
we have an example of a {\em quake-bend} (see Epstein and Marden [\refEpM]). 
We can say that the group $G(\lambda,\tau)$ is obtained from
$G(\lambda,0)$ by doing a quake-bend along $S$ with parameter $\tau$.
That is, for $\lambda\in{\Bbb R}_+$, we take the Fuchsian group 
$G(\lambda,0)$ with generators
$$
S=\pmatrix{ \cosh(\lambda) & \cosh(\lambda)+1 \cr \cosh(\lambda)-1 &
\cosh(\lambda)},\quad
T=\pmatrix{ \coth(\lambda/2) & 0 \cr 0 & \tanh(\lambda/2) }.
$$
This group has a fundamental domain rather like
the one shown in Figure~1.2 except with $\tau=0$ (the copy of the hyperbolic 
plane in question is the hyperplane in ${\Bbb H}^3$ whose boundary is the 
extended real axis). Let $Q(\tau)$ be a loxodromic map with the same fixed
points as $S$ and trace $2\cosh(\tau/2)$. Apply $Q(\tau)$ to that part of
${\Bbb H}^2$ lying above $\ax(S)$, ie those points with 
$\vert z\vert>\coth(\lambda/2)$. What we have done is essentially cut along
$\ax(S)$ and reglued after performing a shear and a bend. Now repeat this
construction along the axis of every conjugate of $S$. This is a quake-bend.
For more details and a precise definition of what is involved, see [\refEpM].
A discussion of quake-bends and
complex Fenchel--Nielsen coordinates in given in Section 5.3 of [\refKSQ].

One can perform this construction for irrational measured laminations.
In this case the new measure is obtained by multiplying the initial bending
measure by the quake-bend parameter. This gives a way of generalising the
Fenchel--Nielsen twist parameter $\tau$ analogous to the way lamination
length generalises the hyperbolic length of a simple closed curve. 

We now relate these ideas by extending the  
{\em $zw=t$--plumbing construction} to this situation. 
Essentially the same construction was used by Earle and Marden [\refEaM] and Kra
[\refKra] in the case of punctured surfaces and it was extended by Ar\'es 
[\refAres] and Parkkonen [\refJouni] for surfaces with elliptic cone points.

Let $X_0$ be a punctured cylinder (as in Section 1). Assume that the boundary 
geodesics $\gamma_1$
and $\gamma_2$ corresponding to boundary components $b_1$ and $b_2$ have
equal length $l=2 \lambda>0$. Let $U_1$ and $U_2$ be neighbourhoods of, 
respectively, the ends of $X_0$ corresponding to $\gamma_1$ and $\gamma_2$. 
Let $\gamma_{12}$ be the
shortest geodesic arc connecting the two boundary components, and let
$$
{\Cal A}_\lambda=\{\zeta\in\Bbb C\mid
e^{-\pi^2/\lambda}<\vert\zeta\vert<1\}
$$
with its hyperbolic metric of constant curvature $-1$. The curve
$\{\vert z\vert=e^{-\pi^2/2\lambda}\}$ is the unique geodesic in ${\Cal
A}_\lambda$ with this metric.

We define local coordinates at the ends of $X_0$ by
$$
z\co U_1\to{\Cal A}_\lambda\quad{\rm and}\quad w\co U_2\to{\Cal A}_\lambda
$$
by requiring that the maps are isometries and that the segments
$\gamma_{12}\cap U_1$
and $\gamma_{12}\cap U_2$ are mapped into ${\Cal A}_\lambda\cap\Bbb
R_+$.
These conditions
define the maps $z$ and $w$ uniquely.

\midinsert
\centerline{
\relabelbox\small
\epsfxsize 9truecm\epsfbox{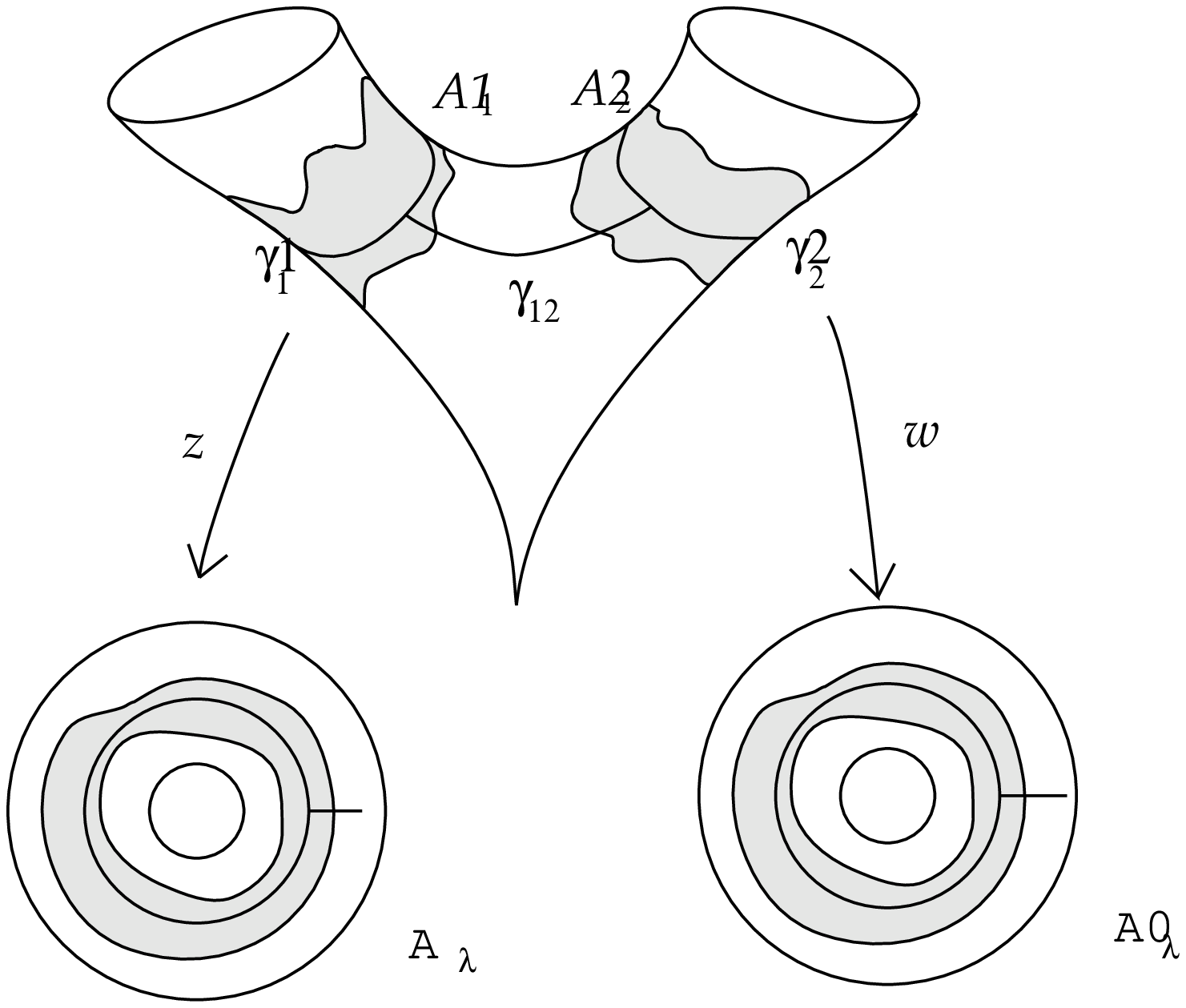}
\relabel {A}{${\Cal A}_\lambda$}
\relabel {A0}{${\Cal A}_\lambda$}
\relabel {A1}{$A_1$}
\relabel {A2}{$A_2$}
\relabel {z}{$z$}
\relabel {w}{$w$}
\relabel {g1}{$\gamma_1$}
\relabel {g2}{$\gamma_2$}
\relabel {g}{$\gamma_{12}$}
\endrelabelbox}
\vglue0.05truein
\centerline{\small Figure~3.1\ssp The $zw=t$ plumbing construction}
\endinsert

If $A\subset X_0$ is an annulus homotopic to a boundary component $b$ of 
$X_0$, we call the component of $\bound A$ separating the other component of 
$\bound A$ from $b$, the outer boundary of $A$. The remaining component of 
$\bound A$ is the inner boundary of $A$.
Assume there are annuli $ A_i\subset U_i$ and a holomorphic homeomorphism
$f\co A_1 \to A_2$ so that
$$
z(x)w(f(x))=t
$$
for some constant
$t\in\Bbb C$ and $f$ maps the outer boundary of $A_1$ to the
inner boundary of $A_2$.
The outer boundaries bound annuli on $X_0$. Remove these annuli to form
a new Riemann surface $X_{{\rm trunc}}$. Define
$$
X_t\co =X_{{\rm trunc}}/\sim,
$$
where the equivalence is defined
by setting 
$$
x\sim y\Longleftrightarrow z(x)w(y)=t.
$$
We say that $X_t$ was obtained from $X_0$ by the $zw=t$
{\em plumbing construction
with plumbing or gluing parameter} $t$. If the annuli $A_i$ can be chosen to be
collar neighbourhoods of the boundary geodesics $\gamma_i$, we say that the
plumbing is {\em tame}.

Next we show that the Fenchel--Nielsen twist parameter is naturally associated 
with a plumbing parameter:

\proclaim{Lemma~3.1}If $G$ is in ${\cal Q}$ with $\lambda\in{\Bbb R}_+$ then
$t=e^{-\pi^2/\lambda}\ e^{-\pi i\tau/\lambda}=e^{i \pi \mu}$ where
$\mu=(i\pi-\tau)/\lambda$.\endproc

\prf Let $\Pi\co {\Bbb H}^2\to {\Bbb H}^2/G_0$ be the canonical projection.
Let $\tilde{\gamma_1}$ be the geodesic in ${\Bbb H}^2$ connecting the fixed points
of $S$ and
$\tilde{\gamma_2}$ the geodesic connecting the fixed points of $S'$ . Now the
boundary
geodesics for which the gluing will be done are
$\gamma_i=\Pi\left(\tilde{\gamma_i}\right)$.
The local coordinates are given by
$$
z(P)= \exp\left({{\pi i}\over \lambda}
\log\left({\Pi^{-1}(P)\sinh(\lambda/2)+\cosh(\lambda/2)}\over
{-\Pi^{-1}(P)\sinh(\lambda/2)+\cosh(\lambda/2)}\right)\right),
$$
and
$$
w(Q)= \exp\left({{\pi i}\over \lambda}
\log\left({\Pi^{-1}(Q)\cosh(\lambda/2)-\sinh(\lambda/2)}
\over
{\Pi^{-1}(Q)\cosh(\lambda/2)+\sinh(\lambda/2)}\right)\right).
$$
Substituting for $T$ we see, after simplifying, that
$$
z\bigl(T(Q)\bigr)= \exp\left({{\pi i}\over \lambda}
\log\left(e^{-\tau}\,{\Pi^{-1}(Q)\cosh(\lambda/2)+\sinh(\lambda/2)}\over
{-\Pi^{-1}(Q)\cosh(\lambda/2)+\sinh(\lambda/2)}\right)\right).
$$
Thus $z(T(Q))\ w(Q)= \exp(-\pi^2/\lambda-\pi i \tau/\lambda)$ as claimed.
\endprf

The same proof also yields the following:

\proclaim{Proposition~3.2}The classical Fenchel--Nielsen construction 
is a $zw=t$ plumbing
construction for a parameter $t$ of modulus $e^{-\pi^2/\lambda}$.
\endproc

\section{$\lambda$--slices}

In this section we keep $\lambda$ real but allow $\tau$ to be complex.
When $\theta=\Im(\tau)$ is in the interval $(0,\pi]$ we will
show that the axis of $S$ is the pleating locus on one component of the
convex hull boundary and when $\theta\in[-\pi,0)$ then it is the pleating 
locus on the other component. We will show that $\tau$ has an interpretation 
as a {\em complex shear} along the pleating locus, $\ax(S)$, see Parker and 
Series [\refPS]. The complex shear $\sigma$ is defined as follows.
The imaginary part of $\sigma$, which we require to be in the interval 
$(-\pi,\pi)$, is the bending angle on the convex hull boundary across 
$\ax(S)$. The real part of
$\sigma$ defined as follows. Let $\eta$ be the unique simple geodesic arc
in the convex hull boundary from $\ax(S)$ to itself and orthogonal to 
$\ax(S)$ at both ends. Then we form a curve in the convex hull boundary in
the homotopy class specified by $T$ by going along $\eta$ and then along
$\ax(S)$. The real part of the complex shear is the signed distance we go 
along $\ax(S)$. This definition is made precise on page 172 of [\refPS].
The theorems of this section should be compared with the 
constructions found in [\refKra] and section 2.2 of [\refKSM]. We also note that
one may use the formulae of [\refPS] to show that, when $\lambda$ is real, the
imaginary part of $\sigma$ cannot be $\pm\pi$, Proposition~7.1 of [\refKSPQ].

\midinsert
\centerline{\leavevmode
\epsfxsize 4 truecm\epsfbox{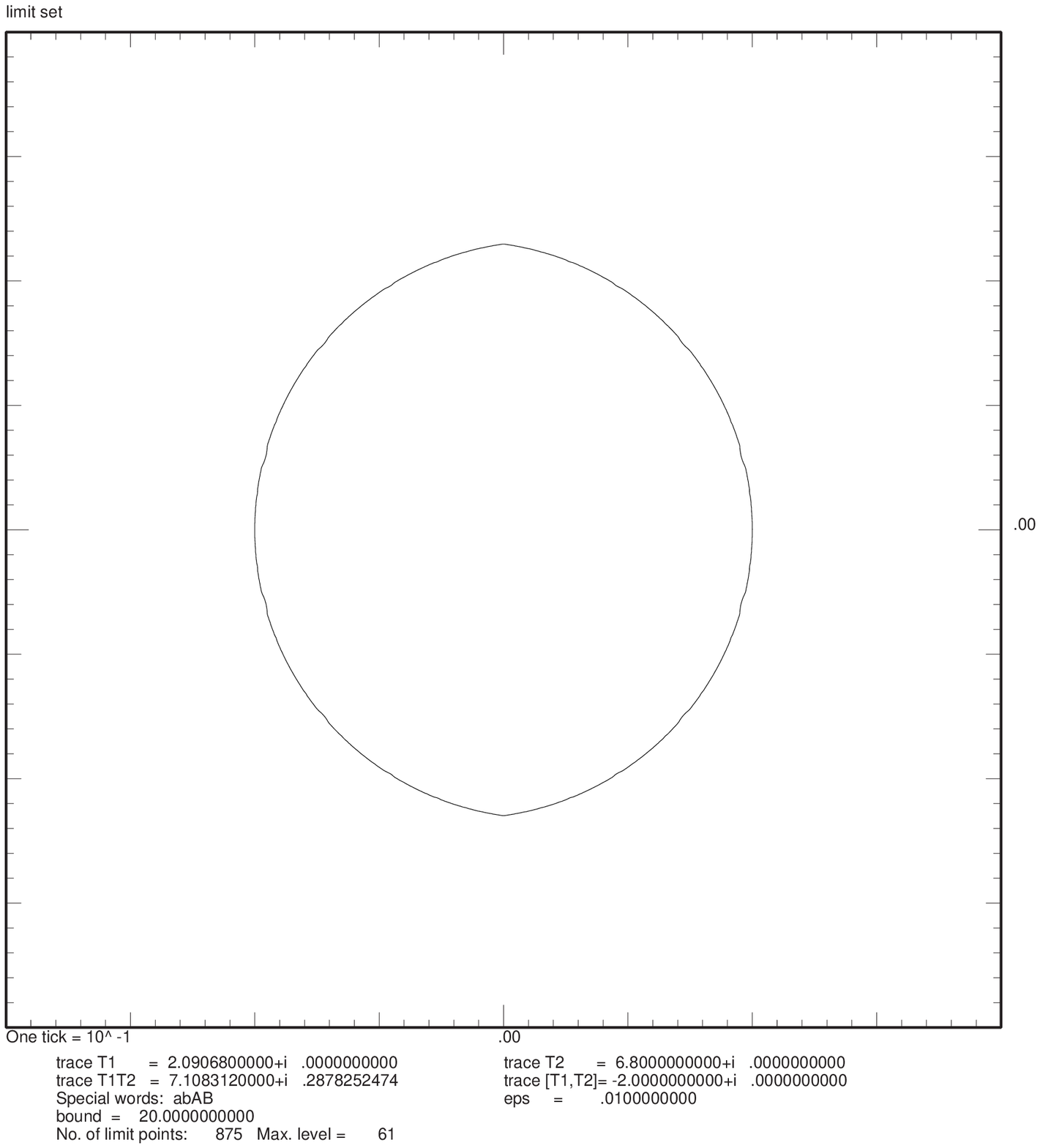}
\epsfxsize 4 truecm\epsfbox{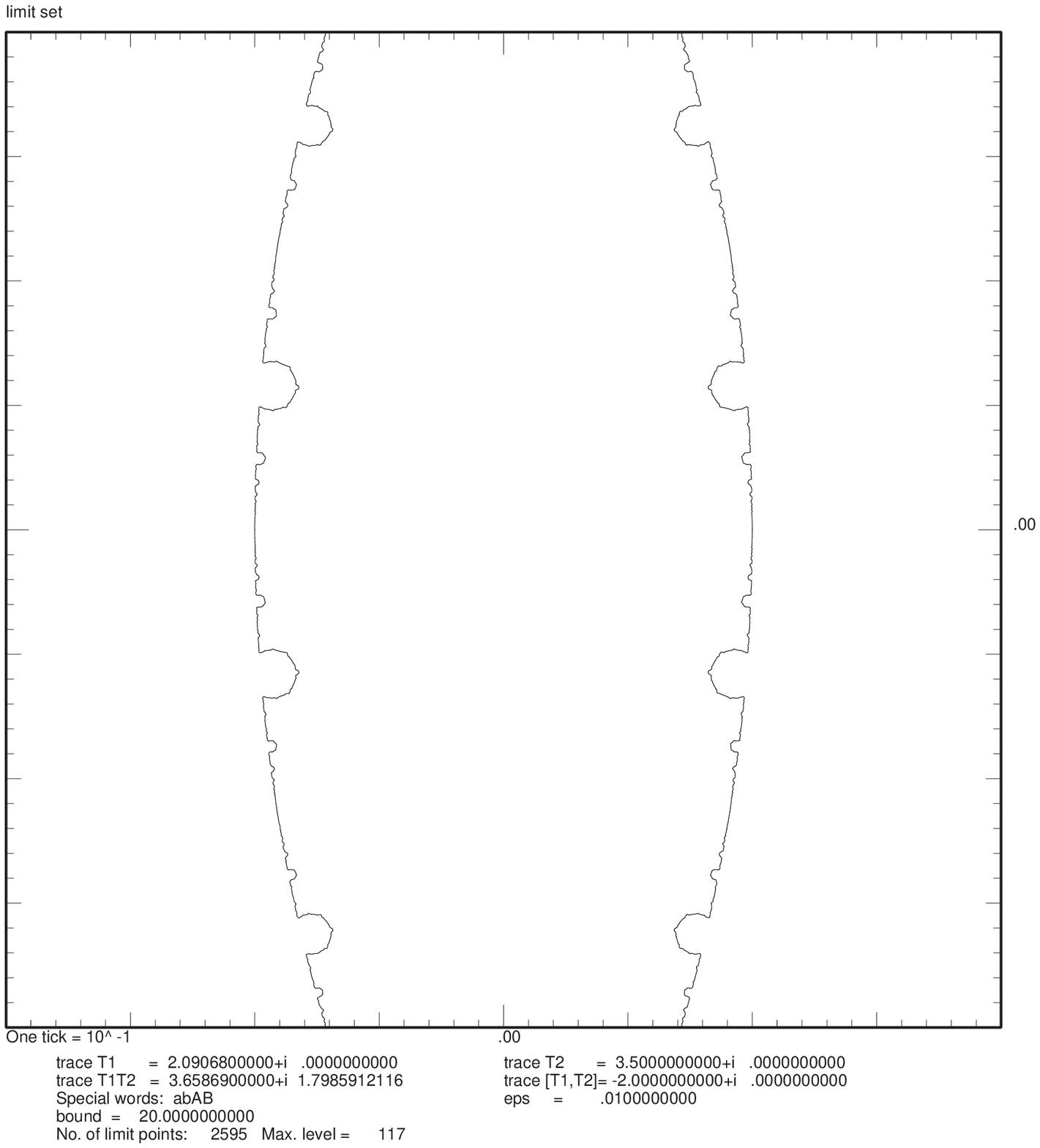}
\epsfxsize 4 truecm\epsfbox{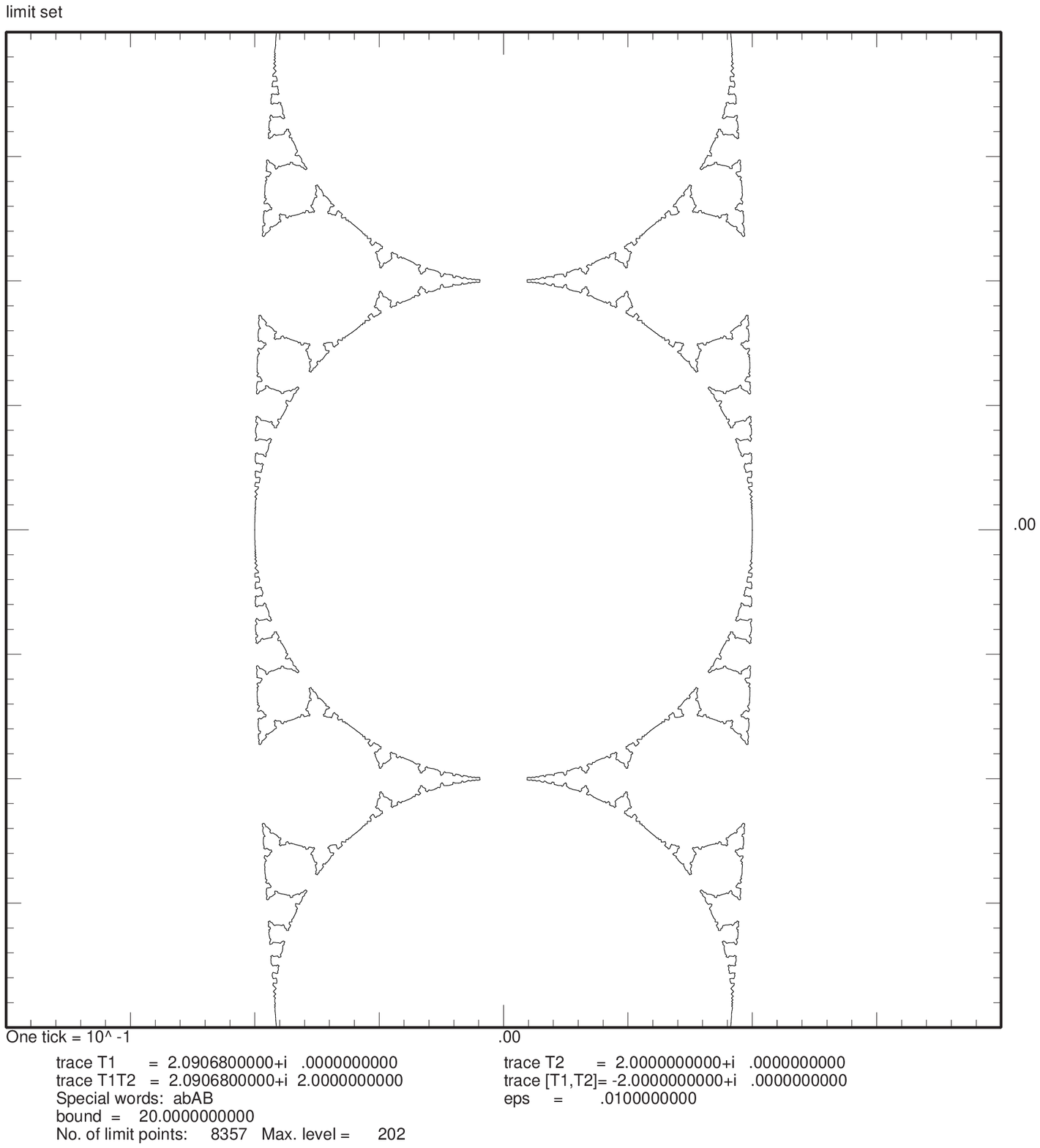} }
\centerline{\small Figure~4.1\ssp Limit sets of groups in a $\lambda$--slice}
\endinsert

Let us fix $\lambda>0$. Consider the set 
$$
\{\tau\in\Bbb C\mid (\lambda,\tau)\in{\Cal FN} \}
$$ 
The {\em $\lambda$--slice} ${\cal Q}_\lambda$ is defined to be the component
of this set containing the points where $\tau\in{\Bbb R}$
(compare with the quake-bend planes of [\refKSPQ]).
We wish to obtain an estimate for the allowed values of $\tau$ for each 
$\lambda$. In order to do this we will construct pleating coordinates on 
each $\lambda$--slice. A first approximation can be achieved by estimating the 
values of $\theta=\Im(\tau)$ that correspond to tame plumbing constructions. 
The following theorem is an explicit version of Theorem~6.1 of [\refKSQ]. 
Specifically, we show that the constant $\epsilon$ of that theorem can be 
taken as $\theta_0=2\arccos\bigl(\tanh(\lambda)\bigr)$ (compare Section 6 of 
[\refCKB]). 
Because the point $(\lambda,i\theta_0)$ is on the boundary of quasi-Fuchsian 
space, there can be no larger uniform bound on $\Im(\tau)$ that ensures 
discreteness. The fact that $\theta$
is the imaginary part of the quake-bend will follow from Theorem~4.2.
 
\proclaim{Theorem~4.1}Let $\theta_0\in(0,\pi)$ be
defined by the equation 
$\cos(\theta_0/2)=\tanh(\lambda)$. Then for 
$\Im(\tau)=\theta\in(-\theta_0,\theta_0)$
the group $G$ is a quasi-Fuchsian punctured torus group.\endproc

\midinsert
\centerline{
\relabelbox
\epsfxsize 9 truecm\epsfbox{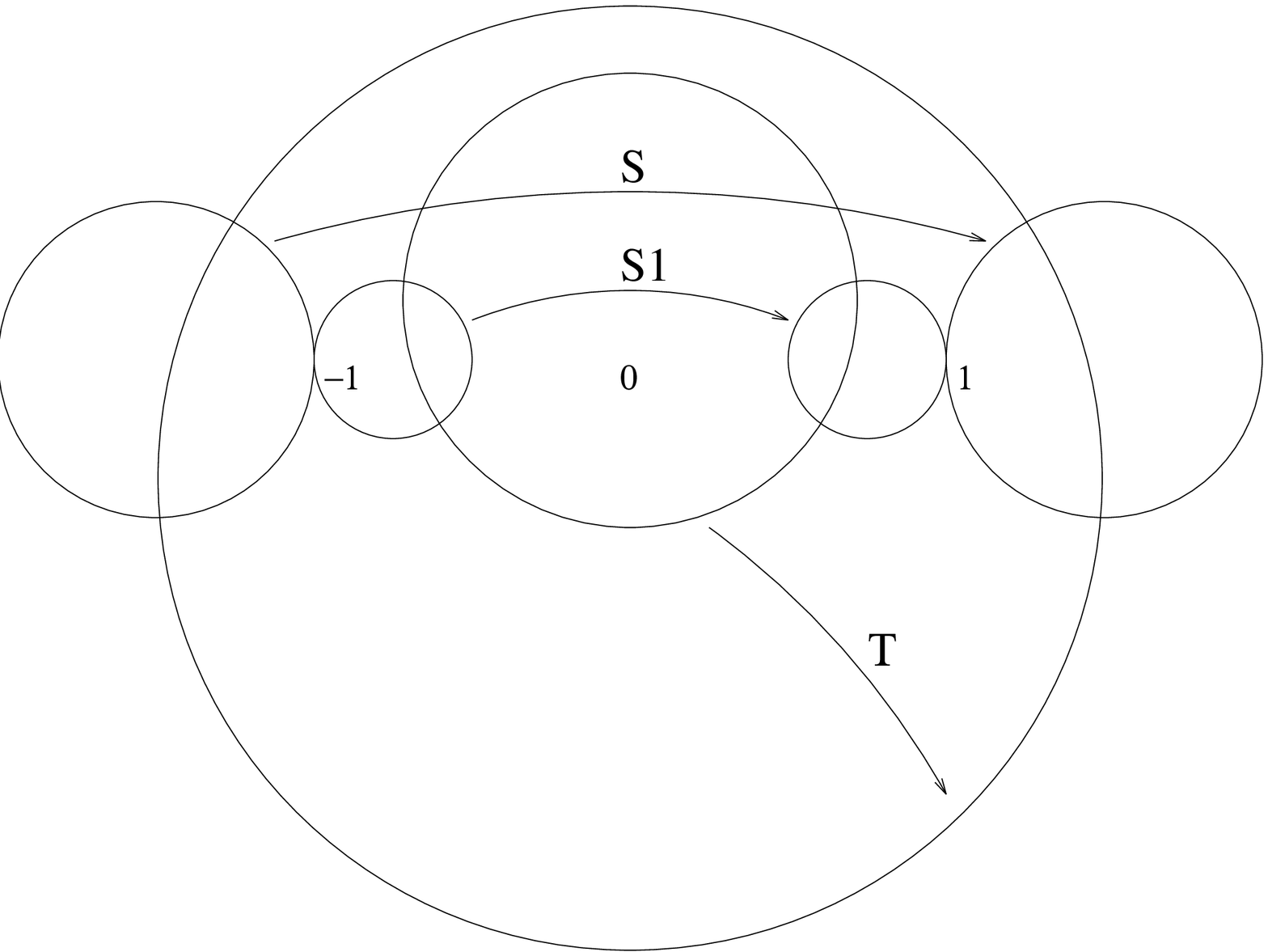}
\relabel{T}{\small$T$}
\relabel{S}{\small$S$}
\relabel{S1}{\small$S'$}
\adjustrelabel <0pt,3pt> {1}{$\scriptstyle1$}
\adjustrelabel <0pt,3pt> {0}{$\scriptstyle0$}
\adjustrelabel <0pt,3pt> {-1}{$\scriptstyle-1$}
\endrelabelbox}
\centerline{\small Figure~4.2\ssp The construction for the combination theorem}
\endinsert

\prf It is easy to check that the circle with centre
at $i\tanh(\lambda/2)\tan(\theta/2)$
and radius $\tanh(\lambda/2)\sec(\theta/2)$ is mapped by $T$ to the circle
with
centre at $-i\coth(\lambda/2)\tan(\theta/2)$ and radius
$\coth(\lambda/2)\sec(\theta/2)$.
Moreover these circles are mapped to themselves under $\gp{S'}$ and $\gp{S}$
respectively (the circles pass through the fixed points of $S'$ and $S$).
Providing the two circles are disjoint then the
annulus between them is a fundamental domain for $\gp{T}$.
It is easy to check that the circles are disjoint if and only if
$\cos(\theta/2)>\tanh(\lambda)$, that is $\theta\in(-\theta_0,\theta_0)$. 
When this happens we can use Maskit's
second combination theorem [\refBMc, \refBM] to show that $G$ is discrete, has a
fundamental domain with two components each of which glues up to give
a punctured torus and $G$ is quasi-Fuchsian.
\endprf

For a positive real number $\lambda$, suppose that $G$ is a 
quasi-Fuchsian punctured torus group. The ordinary set of $G$
has two components. There is an obvious way to label these as the ``top'' and 
``bottom'' components so that, for the case when $G$ is Fuchsian,
the upper half plane is the ``top''component. In what follows, we give a
result that enables us to make this definition precise. Namely in Lemmas~4.3
and~4.4, we show that either the ``top'' component contains the upper half 
plane or the ``bottom''
component contains the lower half plane (or both, in which case the group
would be Fuchsian). When $G$ is strictly quasi-Fuchsian there are two 
components to the convex
hull boundary facing these two components of the ordinary set. We label them
``top'' and ``bottom'' as well (this notation is also used by Keen and Series
on page 370 of [\refKSQ]). Both of these components is a pleated surface and
so we may speak of the pleating locus on the ``top'' and ``bottom''.
The following theorem may be thought of as a generalisation of Proposition~6.2
of [\refPS].

\proclaim{Theorem~4.2}For any parameter in
a $\lambda$ slice ($\lambda\in\Bbb R$) with $\theta\in(0,\pi)$ (respectively
$\theta\in(-\pi,0)$) the pleating locus on the ``bottom'' 
(respectively "top") surface is $S$ and $\tau$ (respectively $-\tau$) is 
the complex shear along $S$ with 
respect to the curve $T$ as defined in [\refPS].\endproc

Intuitively this should be clear as we are keeping $\lambda$ real and
bending away from $\ax(S)$. As we are only bending along one curve the 
result is convex. In the general case we could not expect a Fenchel--Nielsen
complex twist to always be the complex shear on the convex hull boundary
as we may bend along different curves in different directions. In what 
follows we only consider the case $\theta>0$. By symmetry this is sufficient. 
The proof will be by way of several lemmas. 

\proclaim{Lemma~4.3}If $\theta\in(0,\theta_0)$ then the lower half plane
${\Bbb L}$ is contained in $\Omega(G)$.\endproc

\prf
We will consider the lower half plane ${\Bbb L}$ with its Poincar\'e metric. 
We then use plane hyperbolic geometry to prove the result.


Let $D^*$ be the fundamental region for the action of $F=\gp{S,S'}$ on 
${\Bbb L}$ formed by the intersection of ${\Bbb L}$ with the exterior of 
the isometric circles for $S$ and $S'$. That is
$$
D^*=\left\{ z\in {\Bbb L} :
\mod{(\cosh(\lambda)+\eps_1)z+\eps_2\cosh(\lambda)}\ge 1
\hbox{ for all choices of } \eps_1,\eps_2=\pm 1\right\}.
$$
We are now going to consider various hypercycles (that is arcs of circles) 
with endpoints at the fixed points of $S$ and $S'$.
To begin with, let $c_0$ and $c_0'$ be the semicircles centred at $0$ of
radius $\coth(\lambda/2)$ and $\tanh(\lambda/2)$. Clearly these are the
Poincar\'e geodesics joining the fixed points of $S$ and $S'$ respectively.
Let $D_0$ be the subset of ${\Bbb L}$ between these two semi-circles:
$$
D_0=
\left\{z\in {\Bbb L} : \tanh(\lambda/2)\le\mod{z}\le\coth(\lambda/2)\right\}.
$$
The Nielsen region $N(F)$ of $F=\gp{S,S'}$, that is the hyperbolic convex 
hull of $\Lambda(F)$ in ${\Bbb L}$, is 
$$
N(F)=\bigcup_{g\in F}g(D^*\cap D_0).
$$

\midinsert
\centerline{
\relabelbox
\epsfxsize 12.5 truecm\epsfbox{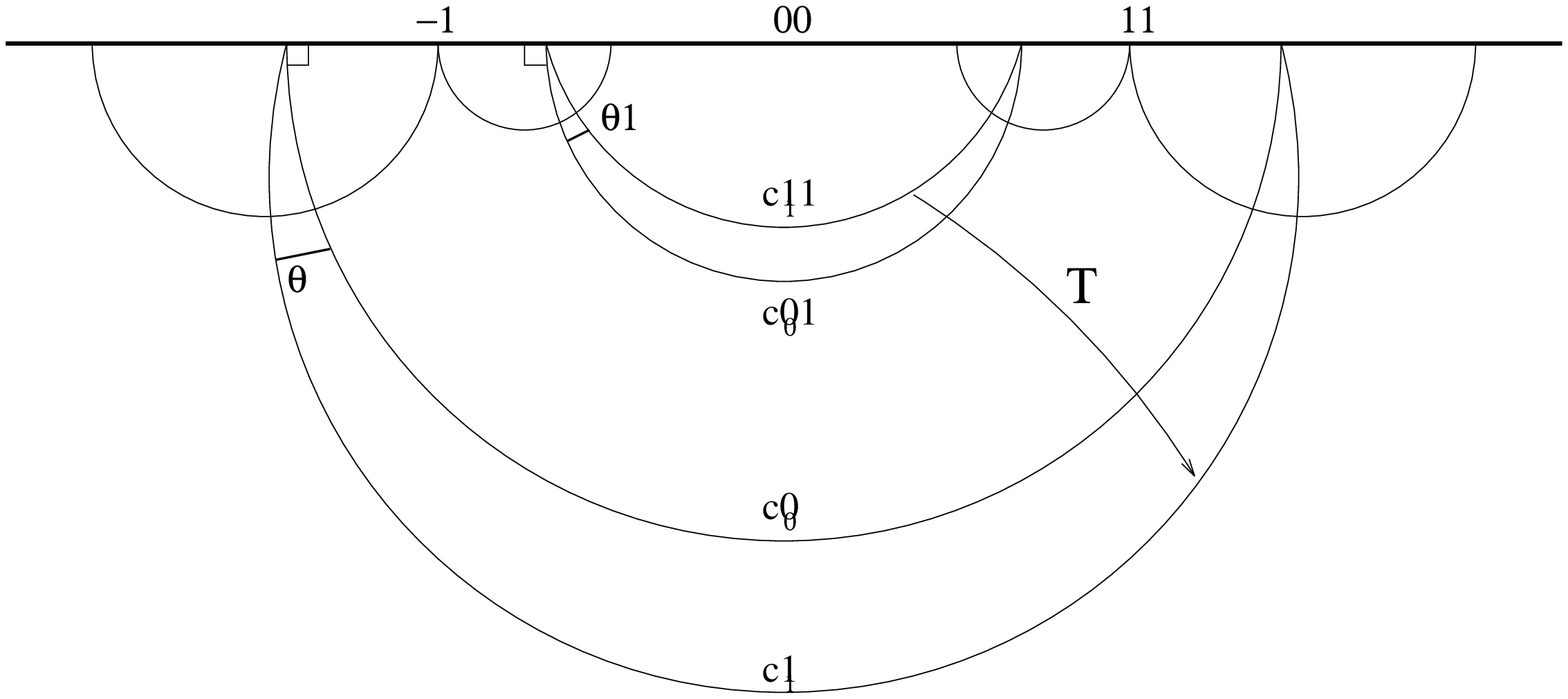}
\relabel {q}{$\theta\over2$}
\relabel {q1}{$\theta\over2$}
\relabel {c0}{\small$c_0$}
\relabel {c1}{\small$c_1$}
\relabel {c01}{\small$c'_0$}
\relabel {c11}{\small$c'_1$}
\relabel {T}{\small$T$}
\relabel {00}{$\scriptstyle0$}
\relabel {11}{$\scriptstyle1$}
\relabel {-1}{$\scriptstyle-1$}
\endrelabelbox}
\vglue 0.05truein
\centerline{\small Figure~4.3\ssp The construction in the lower half plane}
\endinsert

Now consider the circular arcs $c_1$ and $c_1'$ in ${\Bbb L}-D_0$ with 
endpoints at $\pm \coth(\lambda/2)$ and $\pm\tanh(\lambda/2)$ which make an 
angle $\theta/2$ with $c_0$ and $c_0'$ respectively. In other words $c_1$ is
the arc of the circle centred at $-i\coth(\lambda/2)\tan(\theta/2)$ with
radius $\coth(\lambda/2)\sec(\theta/2)$ lying in the lower half plane.
Similarly $c_1'$ is the intersection of ${\Bbb L}$ with the circle centred
at $i\tanh(\lambda/2)\tan(\theta/2)$ with radius 
$\tanh(\lambda/2)\sec(\theta/2)$. Figure~4.3 shows $c_1$ and $c_1'$.
Observe that $c_1$ and $c_1'$ are a constant distance $d(\theta)$ from 
$c_0$ and $c_0'$ where
$$
d(\theta)=\log\bigl(\sec(\theta/2)+\tan(\theta/2)\bigr).
$$
Denote the lune between $c_0$ and $c_1$ by $B(\theta)$ and the lune between
$c_0'$ and $c_1'$ by $B'(\theta)$.
Let $D_1$ be the subset of the lower half plane lying between $c_1$ and $c_1'$.
Now $D_1$ is just the intersection of ${\Bbb L}$ with the fundamental region 
for $T$ considered in Theorem~4.1. One of the consequences of Maskit's 
combination theorem is that $D^*\cap D_1$ is contained in $\Omega(G)$. 
(It is at this point that we have used $\theta<\theta_0$.) Let $N(\theta)$ be
the union of all $F$ translates of $D^*\cap D_1$:
$$
N(\theta)=\bigcup_{g\in F}g(D^*\cap D_1).
$$
It is clear that $N(\theta)$ is just the $d(\theta)$ neighbourhood of $N(F)$.
Since $D^*\cap D_1$ is contained in $\Omega(G)$ then so is $N(\theta)$.

We are going to mimic this construction with more arcs. 
For each $n$ with $n\theta<\pi$, let $c_n$ and $c_n'$
be the circular arcs in ${\Bbb L}-D_0$ with endpoints at $\pm\coth(\lambda/2)$ 
and $\pm\tanh(\lambda/2)$ making an angle of $n\theta/2$ with $c_0$ and
$c_0'$ respectively. That is $c_n$ is the arc of a circle with centre 
at $-i\coth(\lambda/2)\tan(n\theta/2)$ and radius 
$\coth(\lambda/2)\sec(n\theta/2)$ and $c_n'$ is the arc of a circle
with centre at $i\tanh(\lambda/2)\tan(n\theta/2)$ and radius
$\tanh(\lambda/2)\sec(n\theta/2)$.
As before, $c_n$ is a constant distance $d(n\theta)$ from $c_0$ and 
$c_n'$ is the same distance from $c_0'$. We define $D_n$, the subset of
${\Bbb L}$ between $c_n$ and $c_n'$, and the lunes $B(n\theta)$ and
$B'(n\theta)$ as before. Let
$$
N(n\theta)=\bigcup_{g\in F}g(D^*\cap D_n).
$$
Again $N(n\theta)$ is the $d(n\theta)$ neighbourhood of $N(F)$.

Furthermore, let $n_0$ be the integer with $(n_0-1)\theta<\pi\le n_0\theta$. 
We define arcs $c_{n_0}$ and $c_{n_0}'$ which are now in the closed upper 
half plane. We also define $B(n_0\theta)$, $B'(n_0\theta)$ and $N(n_0\theta)$ 
geometrically but remark that these no longer have any metrical properties. 
An important observation is that ${\Bbb L}$ is contained in $N(n_0\theta)$.

\midinsert
\centerline{
\relabelbox
\epsfxsize 12.5 truecm\epsfbox{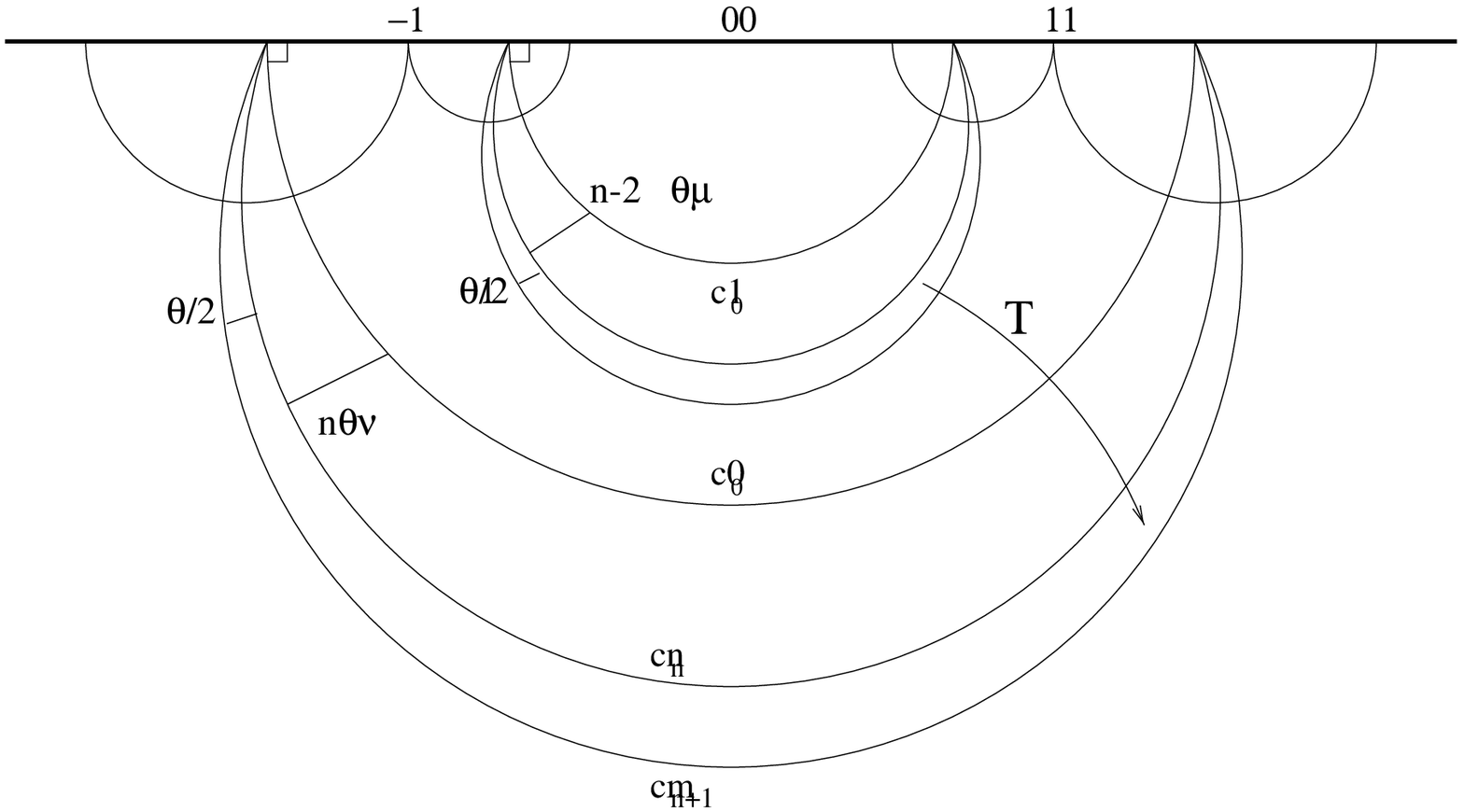}
\relabel {q}{$\theta\over2$}
\relabel {q1}{$\theta\over2$}
\relabel {n}{${\scriptstyle n}{\theta\over2}$}
\relabel {n-2}{${\scriptstyle(n-2)}{\theta\over2}$}
\relabel {cn}{\small$c_n$}
\relabel {cm}{\small$c_{n+1}$}
\relabel {c0}{\small$c_0$}
\relabel {c1}{\small$c'_0$}
\relabel {T}{\small$T$}
\relabel {00}{$\scriptstyle0$}
\relabel {11}{$\scriptstyle1$}
\relabel {-1}{$\scriptstyle-1$}
\endrelabelbox}
\vglue0.05truein
\centerline{\small Figure~4.4\ssp The inductive step}
\endinsert

The rest of the proof follows by an induction from $n=1$ up to
$n=n_0$. We claim that, for $1\le n< n_0$ that if $B(n\theta)$ and 
$B'(n\theta)$ are in $\Omega(G)$ then so are $B\bigl((n+1)\theta\bigr)$ and
$B'\bigl((n+1)\theta\bigr)$. This in turn means that $N\bigl((n+1)\theta\bigr)$
is in $\Omega(G)$. In particular $N(n_0\theta)$, which contains ${\Bbb L}$, 
is in $\Omega(G)$.

Thus all we have do is prove the claim, which we now do. 
Since $B(n\theta)$ and $B'(n\theta)$ are contained
in $\Omega(G)$ then so is $N(n\theta)$. Consider 
$T^{-1}B\bigl((n+1)\theta\bigr)$. Since $c_{n+1}$ makes an angle of 
$n\theta/2$ with $c_1$ and $T$ acts conformally on $\hat{\Bbb C}$ we see that 
$T^{-1}(c_{n+1})$ makes an angle of $n\theta/2$ with $T^{-1}(c_1)=c_1'$,
see Figure~4.4.
In other words $T^{-1}(c_{n+1})$ is a hypercycle a constant distance
$d\bigl((n-1)\theta\bigr)$ from $c_0'$ (also it is not $c_{n-1}'$).
This means that $T^{-1}(c_{n+1})$, and hence also 
$T^{-1}B\bigl((n+1)\theta\bigr)$, is contained within the $d(n\theta)$
neighbourhood of $N(F)$, that is $N(n\theta)$. Since $N(n\theta)$
was assumed to be in $\Omega(G)$, we see that $T^{-1}B\bigl((n+1)\theta\bigr)$
and hence also $B\bigl((n+1)\theta\bigr)$ is contained in $\Omega(G)$, as 
claimed. 
We remark that if $n>n_0$ then $T^{-1}(c_{n+1})$ lies in the closed upper
half plane and the argument breaks down. 
A similar argument shows that $B'\bigl((n+1)\theta\bigr)$ is
also contained in $\Omega(G)$. This completes the proof. \endprf

\proclaim{Lemma~4.4}If $\tau\in{\Cal Q}_\lambda$ and $\theta\in(0,\pi)$ then 
the pleating locus on the ``bottom'' surface is $S$.\endproc

\prf
Suppose first that $\theta\in(0,\theta_0)$. From Lemma~4.3 we see that
${\Bbb L}$ is contained in $\Omega(G)$. Thus the geodesic plane in ${\Bbb H}^3$ 
with boundary the real axis is a support plane for $\partial C(G)$.
Moreover the image of this plane under $T$ must also be a support plane
for $\partial C(G)$. As the intersection of these two planes is the
axis of $S$ we have the result.

Now consider $\tau=t+i\theta\in{\Cal Q}_\lambda$ and 
$\theta\in[\theta_0,\pi)$. 
We proceed as in Proposition~5.4 of [\refKSM]. Suppose that $S$ is not the 
pleating locus for the bottom surface. Consider a path $\alpha$ in
${\Cal Q}_\lambda$ joining
$\tau$ with $\tau'=t'+i\theta'$ where $\theta'\in(0,\theta_0)$. 
Without loss of generality, suppose that if $\tau\in\alpha$ then
$\Im(\tau)\ge \theta'>0$.
We know that at $\tau'$ the pleating locus on the bottom surface is $S$.
Using the standard identification of projective measured laminations
on the punctured torus with the extended real line (with the topology
given by stereographic projection of the usual topology on the circle)
then Keen and Series show that the pleating locus is continuous with
respect to paths in ${\cal Q}$ [\refKSC].
Therefore there are points on the path $\alpha$ for which the pleating locus
is a projective measured lamination arbitrarily close to $\gamma_\infty$.
In particular there are points where the pleating locus is $\gamma_m$
for $m\in{\Bbb Z}$ which corresponds to $W_m=S^{-m}T\in G$ (in the next
section we will give more details of how to associate words with simple
closed curves). In particular, this group element must have real trace.
In other words there is a point of $\alpha$ where 
$\tr(S^{-m}T)=2\cosh(\tau/2+m\lambda)\coth(\lambda)$ is real,
and so
$$
0=\sinh(t/2+m\lambda)\sin(\theta/2).
$$
As $\theta\in[\theta',\pi)$ we see that $\sin(\theta/2)\neq 0$.
Thus $t/2+m\lambda=0$ and $\tr(S^{-m}T)=2\cos(\theta/2)$.
This means $S^{-m}T$ is elliptic and so $\tau$ is not in 
${\Cal Q}_\lambda$ after all.
\endprf

\proclaim{Lemma~4.5}With $S$ and $T$ as in the theorem and $\theta\in(0,\pi)$ 
(respectively $\theta\in(-\pi,0)$) the complex shear $\sigma$ along $S$ with 
respect to $T$ is $\sigma=\tau$ (respectively $\sigma=-\tau$).\endproc

\prf
The trace of $T$ is
$$
\cosh(\tau/2)\bigl(\coth(\lambda/2)+\tanh(\lambda/2)\bigr)
=2\cosh(\tau/2)\coth(\lambda).
$$
Writing $\tr(T)=2\cosh\bigl(\lambda(T)\bigr) $ and
$\tr(S)=2\cosh\bigl(\lambda(S)\bigr)$ the formula $(I)$ of [\refPS] gives
the complex shear along $S$ with respect to $T$ as $\sigma$ where
$$
\eqalign{
\cosh(\sigma/2)&=\cosh\bigl(\lambda(T)\bigr)\tanh\bigl(\lambda(S)\bigr)
\cr
&=\cosh(\tau/2)\coth(\lambda)\tanh(\lambda)\cr
&=\cosh(\tau/2).}
$$
Thus $\sigma$ and $\tau$ agree up to sign and addition of multiples
of $2\pi i$. Since $\Im(\sigma)$ is in $(0,\pi)$ we find that 
$\sigma=\tau$ when $\theta=\Im(\tau)>0$ and $\sigma=-\tau$ when 
$\theta<0$.\endprf

\section{Pleating rays on $\lambda$--slices}

We have shown that on a $\lambda$--slice the pleating locus on one
component of the convex hull boundary is $\gamma_\infty$ which corresponds
to $S$. We now investigate the intersection of each $\lambda$--slice
with the rational pleating plane associated to the simple closed curves
$\gamma_\infty$ and $\gamma_{p/q}$. We call this intersection a 
{\em pleating ray}. Part of the this section will be a justification of
this name. 

In order to obtain pleating rays on each $\lambda$--slice, we follow the 
arguments in [\refKSM],
many of which are inherently two-dimensional in nature. These arguments have
been superseded by more general arguments in [\refKSPQ]. We give these arguments
to help the reader interpret Figure~5.1 and Figure~6.1 without having to 
refer to [\refKSM] or [\refKSPQ]. But, since these arguments are not new, we 
shall not give all the details. Furthermore, we indicate how one may use 
pleating rays on $\lambda$--slices to obtain the rational pleating planes.
This is the simplest part of the construction of pleating coordinates. 
The more complicated parts are treated at length in [\refKSPQ].

\midinsert
\centerline{\epsfxsize 12.5 truecm\epsfbox{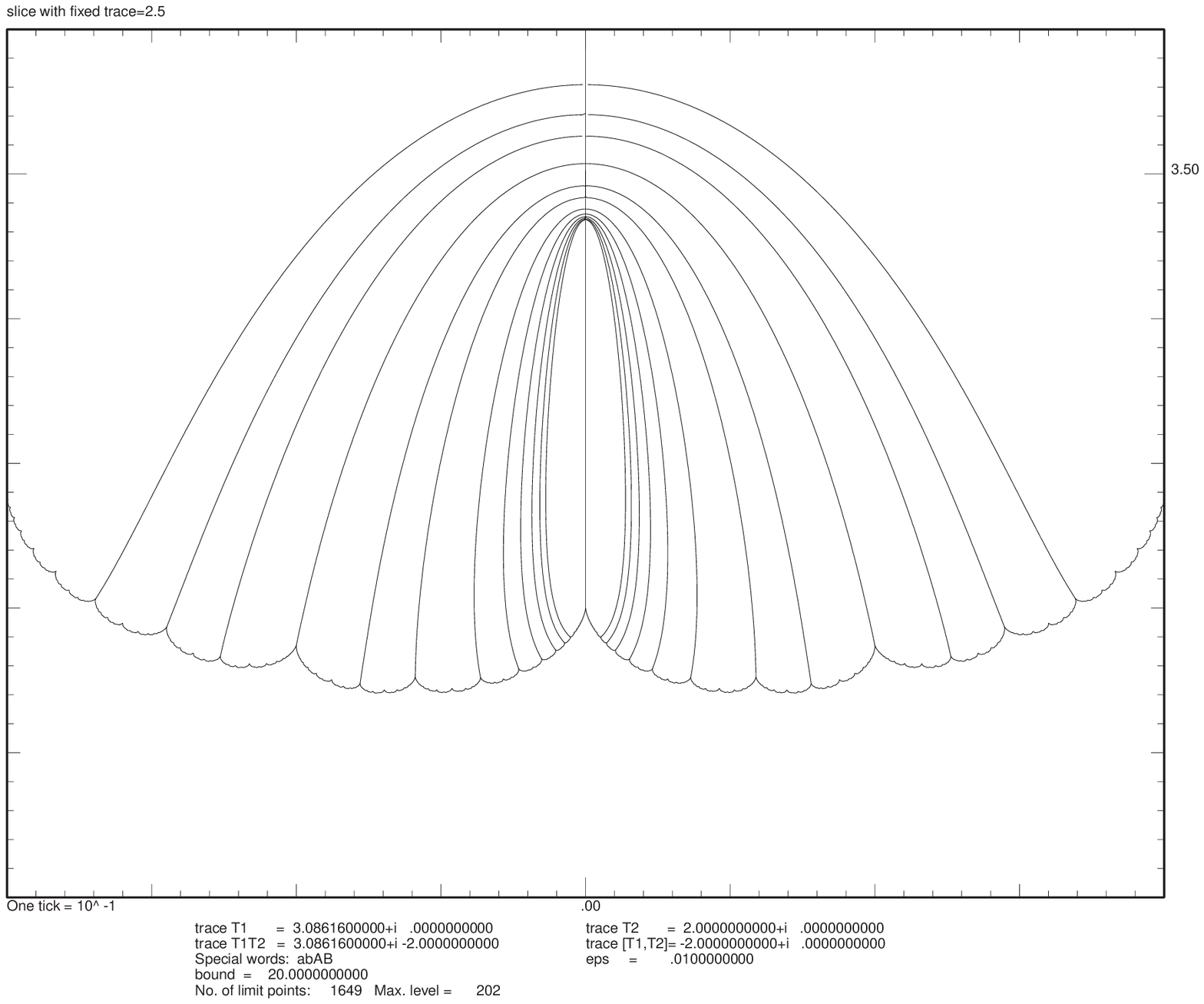} }
\vglue -0.2truein
{\par\leftskip 25pt\rightskip 25pt\small 
Figure~5.1\ssp Part of a slice through ${\Cal Q}$ with $\lambda$ 
held to be real and fixed. In this case $\cosh(\lambda)=5/4$. 
This figure shows the image of the slice under the $2$ to $1$ map
$\tau\longmapsto i\tr T=2i\cosh(\tau/2)\coth(\lambda)=
{10\over 3}i\cosh(\tau/2)$. 
The figure shows pleating rays for this slice, see [\refKSPQ] or Section 6. 
The vertical line from $10i/3$ upwards represents Fuchsian space (which has
been folded onto itself at the point corresponding to a rectangular torus). 
Observe that the pleating rays meet Fuchsian space orthogonally. \par}
\endinsert

In what follows, we assume that the pleating locus on one component of the
convex hull boundary is $\gamma_\infty$, represented by $S$, and the pleating 
locus the other is also a simple closed curve, $\gamma_{p/q}$ for some 
$p/q\in{\Bbb Q}$. There is a special word $W_{p/q}\in G=\langle S,T\rangle$ 
corresponding to the homotopy class of simple closed curves $[\gamma_{p/q}]$.
These words are defined recursively in [\refDW] (see also Section 3.1 of [\refKSM]) 
but of course, we now need to use the generators $S$ and $T$ defined 
(1.1) and (1.3). First, $W_\infty=S^{-1}$, $W_m=S^{-m}T$ for 
$m\in{\Bbb Z}$. If $qr-ps=1$ then
we inductively define $W_{(p+r)/(q+s)}=W_{r/s}W_{p/q}$.

For each $\gamma_{p/q}$ the {\em $p/q$--pleating ray} 
${\cal P}^{\lambda}_{p/q,\infty}$ on ${\Cal Q}_\lambda$ is defined to be the 
those points of ${\Cal Q}_\lambda$ for which
the pleating locus is $\gamma_{p/q}$ on the ``top'' and $\gamma_\infty$
on the ``bottom''. Thus these points have $\Im(\tau)\in(0,\pi)$, Theorem~4.2. 
Likewise ${\cal P}^{\lambda}_{\infty,p/q}$ consists of
those points in ${\cal Q}_\lambda$ where the pleating locus on the ``top''
surface is $\gamma_\infty$ and that on the ``bottom'' is $\gamma_{p/q}$.
Such points have $\Im(\tau)\in(-\pi,0)$.
This discussion may be summarised in the following result which should be
compared to Theorem~5.1 of [\refKSM].

\proclaim{Proposition~5.1}On each $\lambda$--slice ${\cal Q}_\lambda$ and
for $p/q\in{\Bbb Q}$ the pleating rays ${\cal P}^{\lambda}_{p/q,\infty}$ and 
${\cal P}^{\lambda}_{\infty,p/q}$ each consist of a non-empty, connected, non-singular 
arc on which $\tr(W_{p/q})$ is real and which meet
${\cal F}$ orthogonally at the same point from the opposite side. 
Their other end-points lie on the boundary of ${\cal Q}_\lambda$ and at these 
points $\vert\tr(W_{p/q})\vert=2$.\endproc

Some rational pleating rays are shown in the pictures Figures~5.1 and~6.1.
It can be observed that the pleating rays are non-singular connected arcs
that meet Fuchsian space orthogonally.

\medskip{\bf Sketch proof}\stdspace 
This is an adaptation of ideas in [\refKSM] and [\refKSPQ].
First we fix a particular $\lambda$--slice ${\Cal Q}_\lambda$. 
In Theorem~4.2 we showed that $\gamma_\infty$, represented by $S$, is the
pleating locus on one component of the convex hull boundary. For definiteness
we take this to be the ``bottom'' component. By symmetry all our
arguments go through when the pleating loci are the other way round. 

It was shown in Corollary~6.4 of 
[\refPS] that, when the complex shear is purely imaginary, the pleating locus
on the ``top'' component is $T$ (that is $\gamma_0$). Using a change of 
generators (marking) as in Proposition~2.4, it follows that, when the real 
part of the 
complex shear is $-2m\lambda$, for an integer $m$, then the pleating locus on 
the ``top'' component is $S^{-m}T$ (that is $\gamma_m$). Consider the
line where $\Im(\tau)=\theta_0/2$. Such groups are all quasi-Fuchsian
(Theorem~4.1) and at $\tau=-2m\lambda+i\theta_0/2$ the pleating locus is
$\gamma_m$ for $m\in{\Bbb Z}$.
Thus, by the continuity of the pleating locus, see [\refKSC], as we move along
this line we find points whose pleating locus is given by any real parameter.
This shows that any real pleating ray on ${\cal Q}_\lambda$ is non-empty.

It is clear that ${\cal P}^{\lambda}_{p/q,\infty}$ is contained in the real 
locus of $\tr(W_{p/q})$. We now investigate how this real locus
meets Fuchsian space. Any brach of the real locus of $\tr(W_{p/q})$ contained
in ${\cal Q}_\lambda-{\cal F}$ meets ${\cal F}$ in a singularity of 
$\tr(W_{p/q})$. A result of Wolpert, page 226 of [\refSWS], says that the second 
derivative of $\mod{\tr(W_{p/q})}$ with respect to $\tau$ along 
Fuchsian space is strictly positive. (We have used here that $\gamma_{p/q}$ 
and $\gamma_\infty$ are both simple and they intersect.) Thus 
$\tr(W_{p/q})$ has a unique singularity in ${\Cal F}$ and this singularity is
quadratic. Hence the branches of its real locus must meet orthogonally.
In particular there is one branch meeting ${\cal F}$ at this point on
which $\Im(\tau)>0$ and one brach where $\Im(\tau)<0$.

For $0<p/q<1$ the pleating ray ${\cal P}^{\lambda}_{p/q,\infty}$ (which is 
non-empty) must be contained in the
open set bounded by ${\cal F}$, that is $\Im(\tau)=0$; the pleating rays
${\cal P}^{\lambda}_{0,\infty}$, that is $\Re(\tau)=0$, and 
${\cal P}^{\lambda}_{1,\infty}$, that is $\Re(\tau)=-2\lambda$; and the 
boundary of ${\cal Q}_{\lambda}$.
The pleating ray must be a union of connected components of the
intersection of this set with the real locus of $\tr(W_{p/q})$. The
proof of this statement follows Proposition~5.4 of [\refKSM]. A similar argument
has been used in Lemma~4.4 so we will not repeat it. It is clear that if
the pleating locus on the ``top'' is $\gamma_{p/q}$ and if 
$\vert\tr(W_{p/q})\vert>2$ then the group is in the interior of ${\cal Q}$.
Thus, moving along ${\cal P}^{\lambda}_{p/q,\infty}$ in the direction of increasing
$\vert\tr(W_{p/q})\vert$ we cannot reach the boundary of ${\cal Q}$ and so
we must reach ${\cal F}$. It follows that ${\cal P}^{\lambda}_{p/q,\infty}$ is 
connected and non-singular. If not, there would be at least two branches
of ${\cal P}^{\lambda}_{p/q,\infty}$ on which $\vert\tr(W_{p/q})\vert$ is increasing.
But there is only one branch that meets ${\cal F}$, a contradiction.
A similar analysis takes care of other $p/q$.

Finally, when $\vert\tr(W_{p/q})\vert=2$ the pleating ray reaches the
boundary of ${\cal Q}_\lambda$ and the curve $\gamma_{p/q}$ has become 
parabolic. This completes our sketch proof of Proposition~5.1. \endprf

In order to obtain the pleating planes associated to the pairs
$\gamma_\infty$, $\gamma_{p/q}$ we must vary $\lambda$. As we do this,
the pleating rays on each $\lambda$--slice now sweep out the whole pleating 
plane. Keen and Series prove that this gives a connected, non-singular two 
dimensional subset of ${\Cal Q}$. 
In order to obtain pleating planes associated to other pairs of curves we 
use the change of coordinates given in Proposition~2.4. Specifically, if the 
pleating loci we are interested in are $\gamma_{a/b}$ and $\gamma_{c/d}$ which
intersect $q=ad-bc\neq 0$ times then there is a sequence of Nielsen moves
taking the pair $(\gamma_\infty,\gamma_{p/q})$ to the pair 
$(\gamma_{a/b},\gamma_{c/d})$. Associated to these Nielsen moves is
a biholomorphic change of coordinates on ${\Cal Q}$ and the pleating
plane associated to $\gamma_{a/b}$ and $\gamma_{c/d}$ is the image under
this change of coordinates of the pleating plane associated to
$\gamma_\infty$ and $\gamma_{p/q}$.

We conclude this section with a discussion of how one may take data
associated to one component of the convex hull boundary and find information
about the other component. At first sight it does not seem clear how this
could be done. But, at least when the pleating locus on one component
is a simple closed curve, this follows from the relationship between
complex Fenchel--Nielsen coordinates and Keen--Series pleating invariants.
Let $G$ be a strictly quasi-Fuchsian punctured torus group.
Suppose that the pleating locus on one component of the convex hull boundary
is a simple closed curve $\gamma$ of length $\lambda$. Then we can construct 
Fenchel--Nielsen coordinates relative to a generating pair $S$, $T$ where 
$\gamma$ is represented by $S$. The complex Fenchel--Nielsen coordinates 
are given purely in terms of data
associated to the component of the convex hull boundary on which $\gamma$
is the pleating locus. By considering the associated $\lambda$--slice 
${\cal Q}_\lambda$, we can find the Keen--Series pleating invariants for 
$G$ in terms of the complex Fenchel--Nielsen coordinates. We have not
mentioned lamination length on ${\cal Q}_\lambda$ in the above discussion. 
It suffices to remark
that when the lamination on the other component of the convex hull 
boundary is also a simple closed curve given by $W\in G$, then the lamination 
length can be easily found from $\tr(W)$. For irrational pleating rays,
we just use a continuity argument. In particular, we 
can determine information about the pleating on the other component of the
convex hull boundary (this generalises Corollary~6.4 of
[\refPS], where it is shown that if the pleating locus on one component of
$\partial C(G)/G$ is $S$ and the complex shear is purely imaginary then
the pleating locus on the other component is $T$).
Moreover, if the pleating locus on the other component
of the convex hull boundary is also a simple closed curve, we can use a
sequence of Nielsen moves (see Proposition~2.4) to determine the 
Fenchel--Nielsen coordinates with respect to $\delta$. In fact this is
very straightforward.

On the other hand, suppose the pleating locus is an infinite 
measured lamination $\mu$ with lamination length $l_\mu$. The projective class
$(\mu,l_\mu)$ (see [\refKSPQ]) generalises the choice of simple closed curve
with $\delta$--measure  and the hyperbolic length of that curve. It follows 
from the work of Epstein--Marden, [\refEpM], that the group is completely 
determined by $(\mu,l_\mu)$ and the quake-bend parameter $\tau$ (see 
[\refKSQ, \refKSPQ] for a discussion of the quake-bend parameter for 
quasi-Fuchsian punctured torus groups). These 
generalise the Fenchel--Nielsen coordinates for an infinite lamination. 
However, it does not seem that there is a straightforward way to go explicitly
from these parameters to the pleating invariants or to the corresponding 
parameters on the other component of the convex hull boundary.

\section{Degeneration to the Maskit embedding}

In the previous sections we have considered what happens when $\lambda$ is
a fixed real positive number. In this section, we consider what happens
when $\lambda=0$. We should expect the complex shear to tend to $i\pi$ as
$\lambda$ tends to $0$ (compare Theorem~4.1(i) of [\refPS], see Proposition~6.1
below). This means that complex Fenchel--Nielsen coordinates degenerate.
In this section we show that by using the plumbing parameter instead, we 
obtain the {\em Maskit embedding} of Teichm\"uller space, denoted
${\Cal M}$ (see [\refDW, \refKSM]).
This is defined to be the space of free Kleinian groups $G$ on two 
generators $S$, $T$ up to conjugation, such that each group has the following 
properties. First, the generator $S$ and the commutator $K=T^{-1}S^{-1}TS$ 
are both parabolic. Secondly, the components of the ordinary set are of two 
kinds. Namely, a simply connected, $G$--invariant component whose quotient
is a punctured torus; and also infinitely many round discs whose stabilisers 
are thrice punctured sphere groups, all conjugate within $G$. In other words 
these groups are terminal $b$-groups. This space is a holomorphically 
parametrised
copy of the Teichm\"uller space of a punctured torus. There is a standard
normal form for the generators in terms of a parameter $\mu$, see [\refDW, 
\refKSM], which is
$$
S_0=\pmatrix{ 1 & 2 \cr 0 & 1},\quad
T_0=\pmatrix{ -i\mu & -i \cr -i & 0 }.\eqno(6.1)
$$

The goal of this section is to show that as we let $\lambda$ tend
to zero, the normal form for $S$ and $T$ given in (1.1), (1.3) degenerate
to generators of groups in the Maskit embedding (6.1). 
Moreover, the $\lambda$ slices ${\Cal Q}_\lambda$ with their pleating rays
tend to the Maskit embedding with its pleating rays. We illustrate this
with a series of pictures which should be compared to Figure~1 of [\refKSM].
There is a discussion of how the Maskit embedding lies on the boundary
of quasi-Fuchsian space on page 190 of [\refPS].

Consider the limit of $S$ as $\lambda$ tends to zero:
$$
S_0=\lim_{\lambda\rightarrow 0}
\pmatrix{ \cosh(\lambda) & \cosh(\lambda)+1 \cr
\cosh(\lambda)-1 & \cosh(\lambda)}
=\pmatrix{ 1 & 2 \cr 0 & 1}.
$$
Similarly the limit of $S'=T^{-1}ST$ as $\lambda$ tends to zero is:
$$
S'_0=\lim_{\lambda\rightarrow 0}
\pmatrix{ \cosh(\lambda) & \cosh(\lambda)-1 \cr
\cosh(\lambda)+1 & \cosh(\lambda) }
=\pmatrix{ 1 & 0 \cr 2 & 1}.
$$
The parabolic transformations $S_0$ and $S'_0$ generate the level $2$
principal congruence subgroup of $\mathop{\rm PSL}(2,{\Bbb Z})$, 
a torsion-free triangle
group. A comparison of the plumbing parameter calculated in Lemma~3.1 with 
the corresponding result for terminal $b$-groups (see Kra [\refKra; Section 6.4]) 
suggests that, in order to study the degeneration of quasi-Fuchsian groups in
$\bigcup_{\lambda>0}\Cal Q_\lambda$ as
$\lambda\longrightarrow 0$, it is useful to make a change of parameters
$$
\mu={i\pi-\tau\over\lambda}.
$$
We refer to $\mu$ as the {\em plumbing parameter}. In terms of this parameter 
the matrix $T$ can be written as
$$
T=\pmatrix{ -i\sinh(\lambda\mu/2)\coth(\lambda/2) &
-i\cosh(\lambda\mu/2) \cr
-i\cosh(\lambda\mu/2) & -i\sinh(\lambda\mu/2)\tanh(\lambda/2)}.\eqno(6.2)
$$
Using Lemma~3.1, we see that (1.1) and (6.2) give a parametrisation of 
the generators of $G$ in terms of a length parameter and a plumbing parameter.
The following result on the limit groups, which should
be compared to Theorem~4.1(i) of [\refPS], now follows rather easily:

\proclaim{Proposition~6.1}Consider a sequence of groups where $\lambda$ tends 
to
zero but $\mu$ remains fixed. Then the complex shear along $S$ tends
to $i\pi$.\endproc

\prf The conclusion is immediate from the definition of $\mu$:
$\tau=i\pi-\mu\lambda\to i\pi$ as $\lambda\to 0$.
\endprf

We now show that when $\lambda$ tends to zero with $\mu$ being kept fixed 
we obtain the standard form for group generators in the Maskit embedding.

\proclaim{Proposition~6.2}Assume that $\mu\in\Cal Q_\lambda$ for small
$\lambda$. As $\lambda$ tends to zero the group with parameter
$(\lambda,\mu)$ tends to the terminal $b$-group representing punctured torus on
its invariant component with parameter $\mu$. \endproc

\prf We have already seen that $S_0$ and $S'_0$ have the correct form.

Let $\mu$ be fixed. For small $\lambda$ we have
$$
\sinh(\lambda\mu/2)\coth(\lambda/2)=\bigl(\lambda\mu/2+O(\lambda^2)
\bigr)\bigl(2/\lambda+O(1)\bigr)=\mu+O(\lambda).
$$
Therefore we have
$$
\lim_{\lambda\rightarrow
0}\bigl(\sinh(\lambda\mu/2)\coth(\lambda/2)\bigr)=\mu.
$$
This means that the limit as $\lambda$ tends to zero of $T$ is
$$
\eqalign{
T_0&=\lim_{\lambda\rightarrow 0}
\pmatrix{ -i\sinh(\lambda\mu/2)\coth(\lambda/2) &
-i\cosh(\lambda\mu/2) \cr
-i\cosh(\lambda\mu/2) & -i\sinh(\lambda\mu/2)\tanh(\lambda/2)}\cr
&=\pmatrix{ -i\mu & -i \cr -i & 0 }.\cr}
$$
The limiting matrices $S_0$ and $T_0$ are just the usual group generators 
of terminal $b$-groups in the
Maskit embedding $\Cal M$ of Teichm\"uller space of the punctured torus.
\endprf

The convergence of $\lambda$--slices to $\Cal M$ is illustrated in Figure~6.1.

\midinsert
\centerline{\leavevmode
\epsfxsize 6 cm\epsfbox{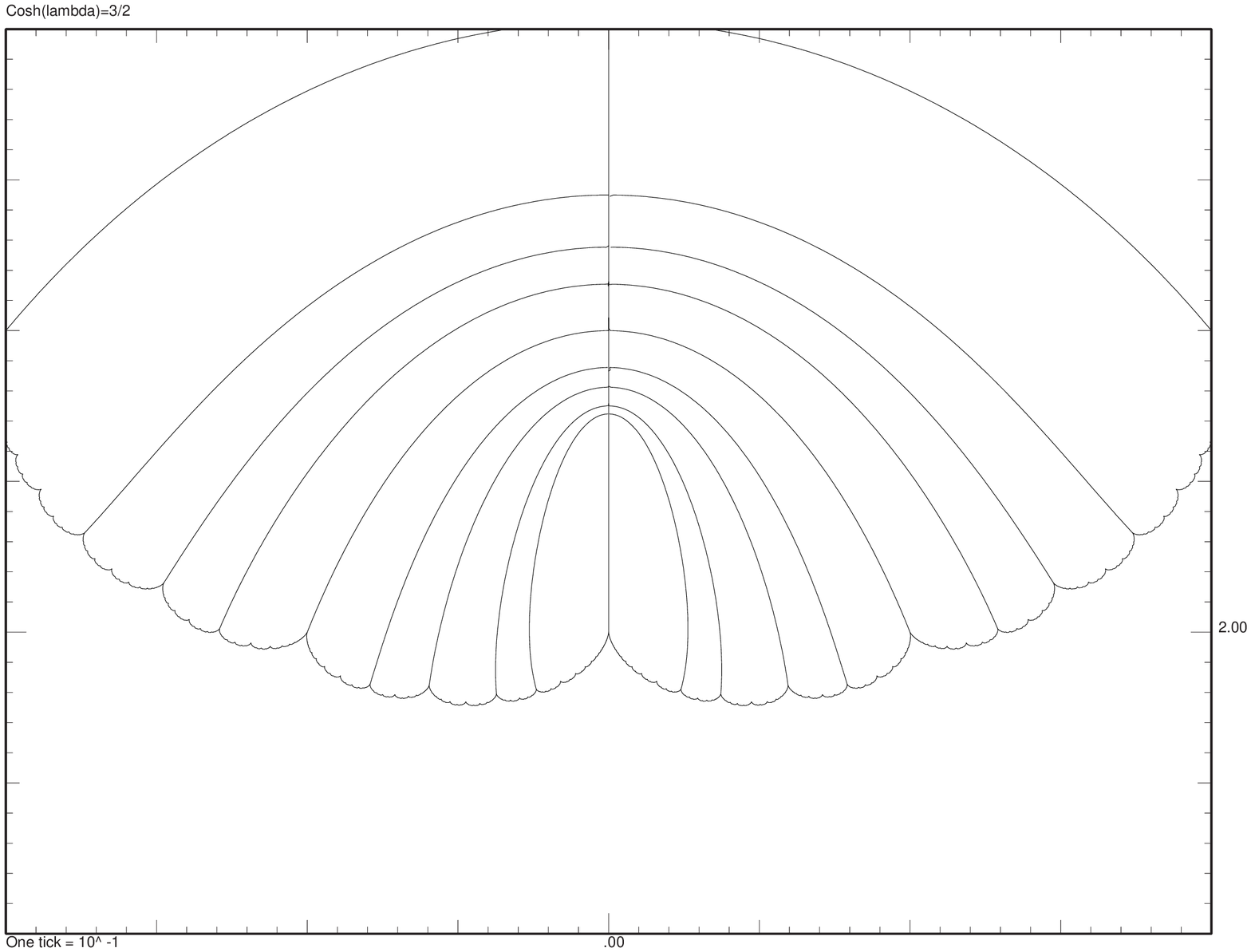}
\epsfxsize 6 cm\epsfbox{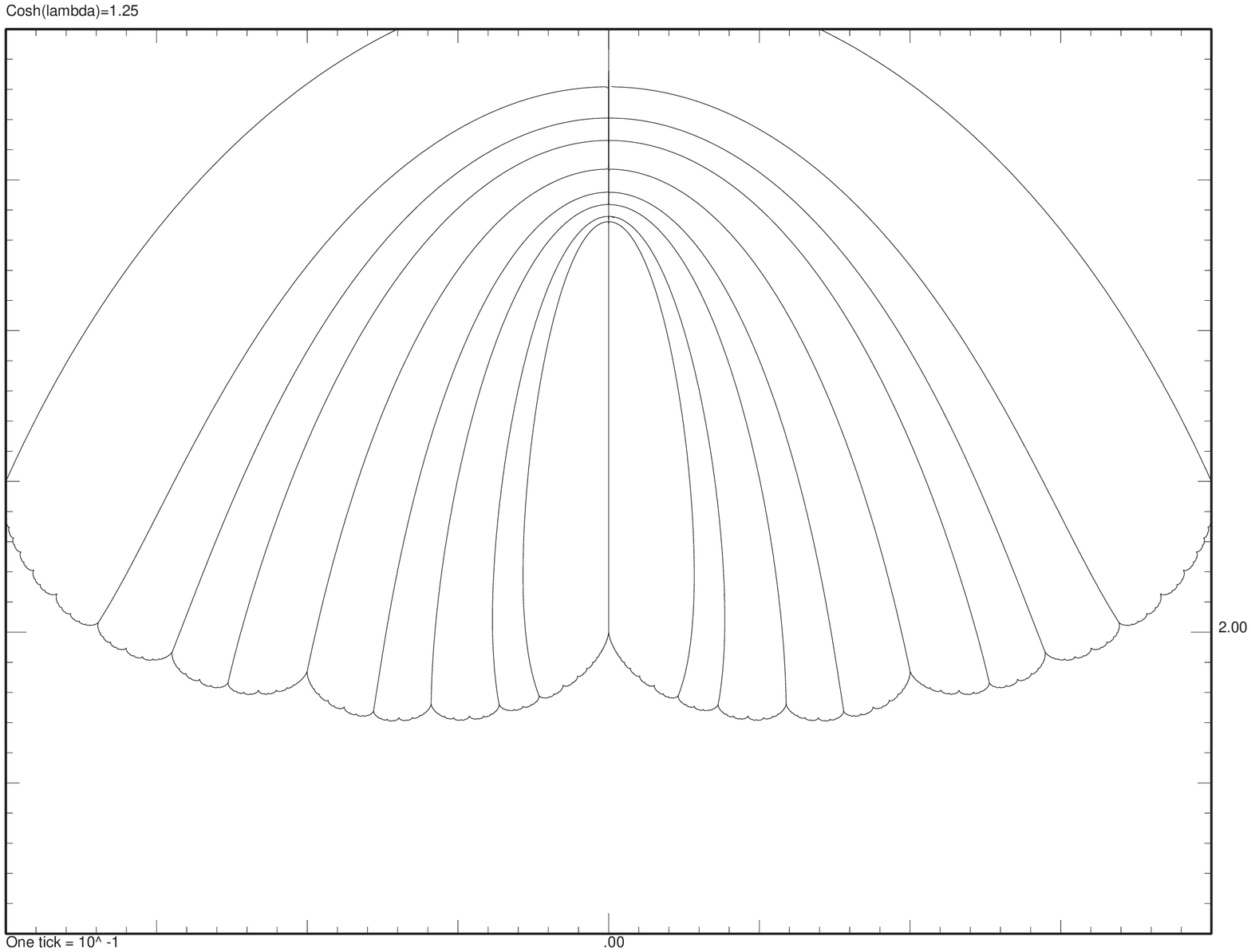} }
\centerline{\leavevmode
\epsfxsize 6 cm\epsfbox{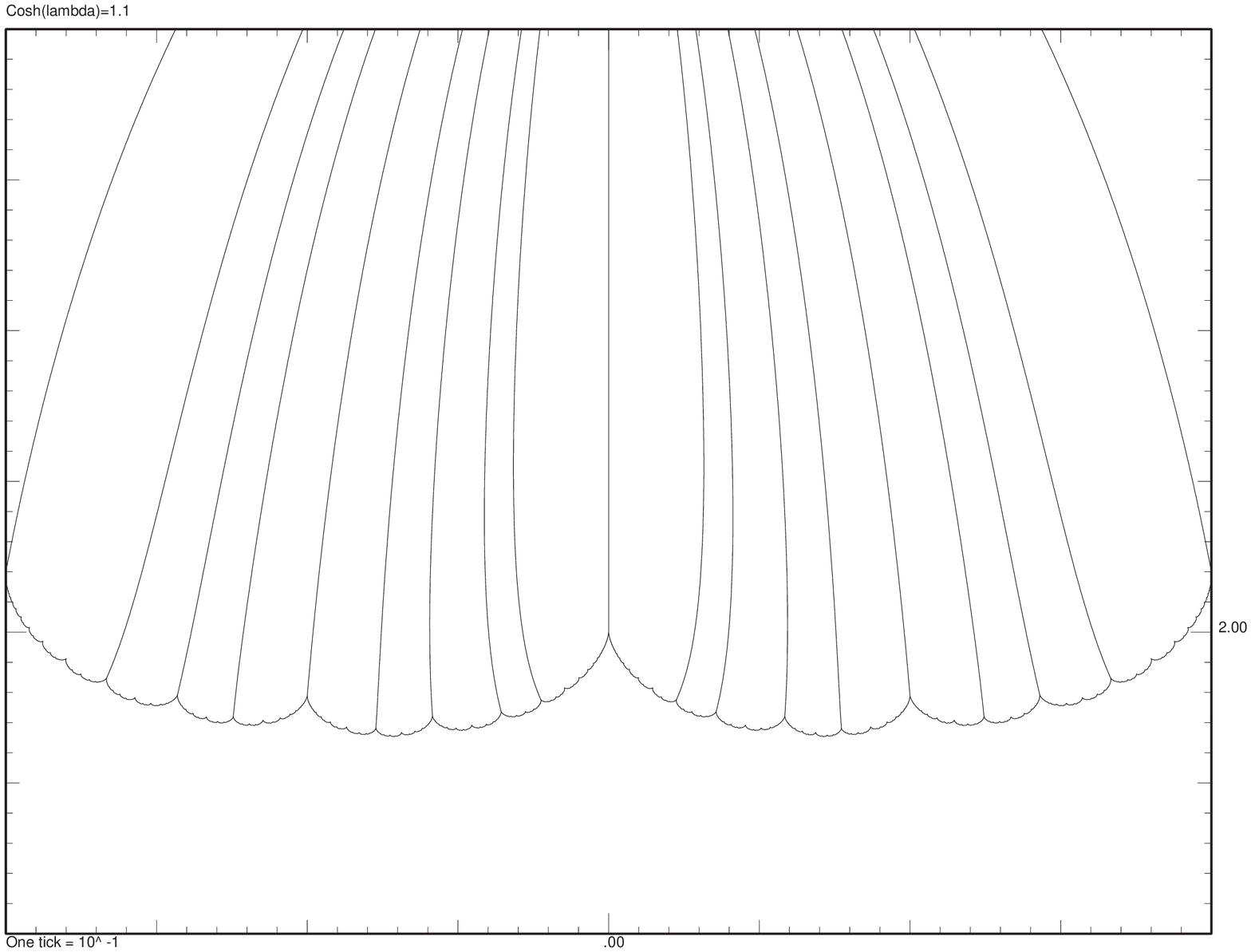}
\epsfxsize 6 cm\epsfbox{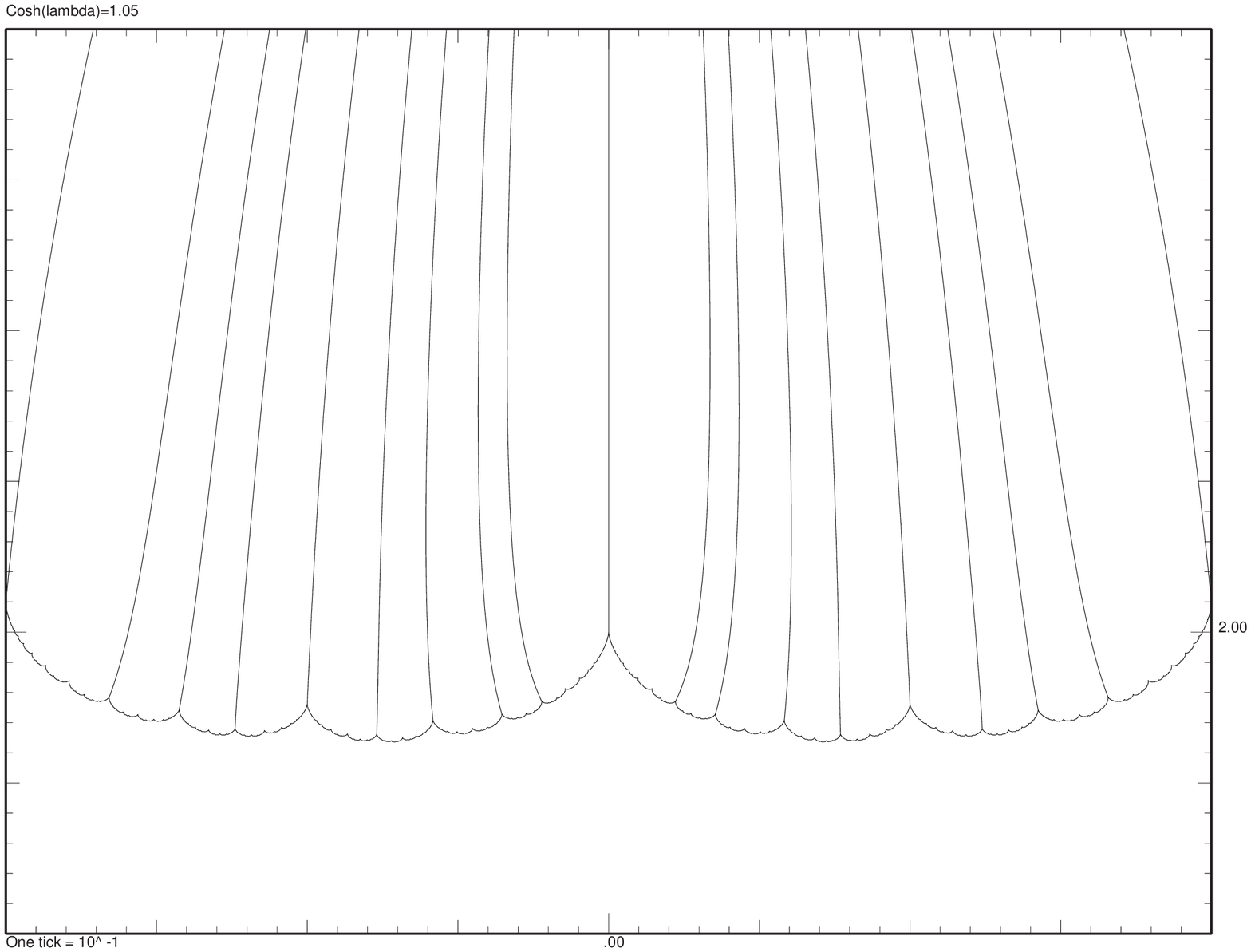} }
{\par\leftskip25pt\rightskip25pt\small
Figure~6.1\ssp $\lambda$--slices for 
$\tr(T)=3$,~$2.5$,~$2.2$ and~$2.1$ drawn with a collection of rational 
pleating rays\par}
\endinsert

\noindent
{\bf Remarks~6.3}\stdspace (a)\stdspace The plumbing construction is
tame when $\Im(\tau)=\theta\in(0,\theta_0)$ or equivalently
$\Im(\mu)\in\bigl((\pi-\theta_0)/\lambda,\pi/\lambda\bigr)$. For small
$\lambda$ we have $\theta_0=\pi-2\lambda+O(\lambda)^2$.  As $\lambda$
tends to zero this interval tends to $(2,\infty)$, which is the
condition for tame plumbing in the Maskit slice, Section 6.2 of
[\refKra] or Proposition~2.3 of [\refDW].

(b)\stdspace  In the $(\lambda,\mu)$ parameters, Fuchsian space corresponds to
the union of the lines $\Im(\mu)=\pi/\lambda$. When $\lambda\to 0$, $\Im\mu\to
\infty$, that is, the closure of Fuchsian space touches $\Cal M$ at the
boundary point corresponding to the parameter $\mu=\infty$ 
(see page 191 of [\refPS]).

\references

\Addresses\recd
\bye